\documentclass[12pt,a4paper,twoside]{article}
\usepackage[latin1]{inputenc}
\usepackage[T1]{fontenc}
\usepackage{graphicx}
\usepackage{pstricks}
\usepackage{pstricks-add}
\NormalCoor
\usepackage{float}

\usepackage{subfigure}
\usepackage{latexsym}
\usepackage{amsmath,amssymb,amsthm}
\usepackage{bbm}
\usepackage{ifpdf}
\usepackage{exscale}
\usepackage{fancyhdr}

\usepackage[DIV12]{typearea}

\newcommand{\dt}{\ensuremath{\textup{d}}}

\newtheorem{theorem}{Theorem}[section]
\newtheorem{lemma}{Lemma}[section]
\newtheorem{proposition}{Proposition}[section]
\newtheorem{corollary}{Corollary}[section]
\newenvironment{prooof}{\noindent{\bf Proof:}}{%
  \hspace*{\fill}$\Box$\par\vskip2ex}
\newenvironment{prooof2}{\noindent{\bf Proof }}{%
  \hspace*{\fill}$\Box$\par\vskip2ex}

\newcommand{\pP}{\ensuremath{\mathbb{P}}}
\newcommand{\pQ}{\ensuremath{\mathbb{Q}}}
\newcommand{\pE}{\ensuremath{\mathbb{E}}}
\newcommand{\pF}{\ensuremath{\mathcal{F}}}

\newcommand{\vp}{\ensuremath{\varphi}}
\newcommand{\ump}{\ensuremath{\mathcal{P}^\downarrow}}
\newcommand{\massp}{\ensuremath{\mathbb{S}_{\leq 1}}}
\newcommand{\omp}{\ensuremath{\mathcal{P}^\circ}}
\newcommand{\osec}{\ensuremath{\mathcal{S}^\downarrow}}
\newcommand{\cs}{\ensuremath{\mathcal{C}}}

\newcommand{\eql}{\ensuremath{\overset{\mbox{\scriptsize\textup{d}}}{=}}}

\theoremstyle{definition}
\newtheorem{remark}{Remark}[section]
\newtheorem{definition}{Definition}[section]

\renewcommand{\theenumi}{{\rm (\roman{enumi})}}

\begin{document}
\pagestyle{fancy} \renewcommand{\headrulewidth}{0mm}
\chead[\leftmark]{\rightmark} 
\rhead[]{\thepage} \lhead[\thepage]{} \cfoot{} 

\pagenumbering{arabic} 
\setcounter{page}{1}
\title{On a ternary coalescent process} \author{Erich Baur\thanks{Institut
    f\"ur Mathematik, Universit\"at Z\"urich, Winterthurerstrasse 190,
    CH-8057 Z\"urich, Switzerland. Email: erich.baur@math.uzh.ch.}}  \date{}
\maketitle \thispagestyle{empty}
\begin{abstract}
  We present a coalescent process where three particles merge at each
  coagulation step. Using a random walk representation, we prove duality
  with a fragmentation process, whose fragmentation law we specify
  explicitly. Furthermore, we give a second construction of the coalescent
  in terms of random binary forests and study asymptotic
  properties. Starting from $N$ particles of unit mass, we obtain under an
  appropriate rescaling when $N$ tends to infinity a well-known binary
  coalescent, the
  so-called  standard additive coalescent.
\end{abstract}

\section{Introduction}
Generally speaking, a stochastic coalescent is a Markov process
describing the coagulation of particles characterized by their size
only. The rate at which particles merge depends just on the members
involved. Conversely, fragmentation processes describe a Markovian
evolution of particles which split independently into new particles
(branching property).  The goal of this paper is to study the stochastic
coalescent with ternary coagulation kernel
\begin{equation*}
\kappa(r,s,t) =r+s+t+3,\quad r,s,t > 0,
\end{equation*}
to which we will simply refer to as ternary coalescent or ternary
coalescent process. Here, three particles of sizes (masses) $r,s,t$ coagulate
into a new particle of size $r+s+t$ at rate $r+s+t+3$.  Although at first
glance, the kernel $\kappa$ may look somewhat arbitrary (for example, it
is not scale invariant), the corresponding
process enjoys rather interesting properties. Similar to the additive
coalescent, that is the coalescent where two particles with masses $s$,
$t$ merge at rate $\tilde{\kappa}(s,t)=s+t$, the state chain of the ternary
coalescent admits different representations. In the spirit of
Bertoin~\cite{BER:4}, we show how it can be obtained by looking at
excursion intervals of a one-dimensional conditioned random walk. As a
by-product of our representation, we establish duality with a 
fragmentation process via time-reversal. We stress that this is a unusual
feature, because the branching property normally fails when time is
reversed in a coalescent process. Section 7 of Bertoin~\cite{BER:3} gives a
brief overview over cases where such a duality relation has been
proven. See also Chapter 5.5 in Pitman's lecture notes~\cite{PIT:StF} for
further discussions.

Using the same construction, we study asymptotic properties of the
ternary coalescent starting from $N$ particles of unit mass. Properly
rescaled in space and time, we observe in the limit $N\rightarrow\infty$
the so-called standard additive coalescent, which has been obtained by
Evans and Pitman in~\cite{EVPIT} as the weak limit $n\rightarrow\infty$ of the
(binary) additive coalescent, started at time $-(1/2)\ln n$ with $n$ atoms of 
size $1/n$. Here, Bertoin's characterization~\cite{BER:2} of the dual
fragmentation process connected to the standard additive coalescent by
time-reversal plays a pivotal role. We emphasize that even though $\kappa$ is
a ternary coagulation kernel, we end up in the limit with a binary
coagulation process. 

We also highlight a second construction of the ternary coalescent
involving random binary forests, following the ideas of Pitman
in~\cite{PIT:99}. In a final remark, we point out that this representation
could instead be used to work out our results. Moreover, we outline a
possible extension of the results to certain $k$-ary coalescent processes.

The rest of this paper is organized as follows. In the first section,
we describe the semigroup of the ternary coalescent and derive some
further properties. We finish this part by computing the one-dimensional statistics for the
underlying state chain starting from an odd number of particles of unit
mass. Its special form already hints at a connection to hitting times of
a one-dimensional nearest neighbor random walk, which we elaborate in the next
section. There we prove duality via time-reversal with a
fragmentation process, using an explicit construction of the coalescent in
terms of ladder epochs. In the third part, we turn our 
attention to random binary trees and find a second interpretation of the
ternary coalescent which is based on random binary forests. Finally we use
again the random walk representation to study asymptotic properties of the coalescent
in the last section.

\section{Some basic properties}
\label{tc}
Throughout this text, let
\begin{equation*}
  \mathbb{N} = \{1,2,\ldots\},\quad\mathbb{Z} = \{\ldots,-1,0,1,\ldots\},\quad\mathbb{Z}_+ = \mathbb{N}\cup\{0\}.
\end{equation*}
The coalescent process will take values in the space of decreasing
numerical sequences with finitely many non-zero terms
\begin{equation*}
 \osec= \left\{{\bf s}=(s_1,s_2,\ldots) : s_1\geq
    s_2\geq\ldots\geq 0, s_k = 0\mbox{ for } k\mbox{ sufficiently large}\right\}.
\end{equation*}
We may think of elements of a sequence ${\bf s}\in\osec$ as (sizes of)
atoms or particles and simply identify ${\bf s}$ with its non-zero
components. If we write ${\bf s} = (s_1,\ldots,s_l)$, the non-zero
components of ${\bf s}$ are precisely given by $s_1,\ldots,s_l$. If ${\bf
  s}=(s_1,s_2,\ldots)\in\osec$ and $1\leq i<j<k$, 
we use the notation ${\bf s}^{i\oplus j\oplus k}$ for the sequence in
$\osec$ obtained from ${\bf s}$ by merging its $i$th, $j$th and $k$th
terms, that is one removes $s_i$, $s_j$, $s_k$ and rearranges the remaining
elements together with the sum $s_i+s_j+s_k$ in decreasing order.

Let us define the object of our interest. Recall the kernel $\kappa$ from
the introduction.
\begin{definition}
  The ternary coalescent with values in $\osec$ and kernel $\kappa$ is a
  continuous time Markov process $\mathcal{X}=(\mathcal{X}(t),\, t\geq 0)$
  with state space $\osec_{'}$ for an appropriate subset $\osec_{'}$
  of $\osec$, and jump rates
\begin{equation*}
  q({\bf s},\cdot) = \sum_{1\leq i < j <k,\,
    s_k>0}\kappa(s_i,s_j,s_k)\delta_{{\bf s}^{i\oplus j\oplus k}}.
\end{equation*}
\end{definition}
This definition can be adapted in an obvious way to other coagulation kernels,
leading to different stochastic coalescent models, for example the additive
coalescent with kernel $\tilde{\kappa}(s,t) = s+t$.

Before looking at concrete realizations, we collect in this section some
basic properties which can be read off from the kernel $\kappa$ and the
very definition of jump-hold processes of the above type. Denote by
$\mathcal{X}=(\mathcal{X}(t),\, t\geq 0)$ the ternary coalescent, started
from a finite configuration ${\bf r}=(r_1,\ldots,r_N)\in\osec$, where
$N=2n+1$, $n\in\mathbb{Z}_+$. We write $M=r_1+\ldots+r_N$ for the total
mass in the system. For every $k=0,\ldots,n+1$, let $T_k$ be the instant of
the $k$th coagulation, with the convention $T_0=0$, $T_{n+1} = \infty$. The
state chain or skeleton chain $\mathcal{X}'$ of the coalescent process is
given by $\mathcal{X}_k'=\mathcal{X}(T_k)$, $k=0,\ldots,n$. We use the
expression $\#(t)$ for the number of particles at time $t$, whereas $J(t)=
\max\{k\in\mathbb{Z}_+ : T_k\leq t\}$ stands for the number of jumps up to
time $t$. Note that $\#(t) = N-2J(t)$.
\subsection{State chain and semigroup}
\begin{proposition}
\label{tc-generalprop}
In the preceding notation, the following holds true.
\begin{enumerate}
 \item The sequence $\Delta_k = T_k-T_{k-1}$, $k=1,\ldots,n$, of the waiting
  times between two coagulations is a sequence of independent exponential
  variables with respective parameters 
\begin{equation*}
\alpha(k)=\frac{1}{2}(M+N+2-2k)(N+1-2k)(N-2k).
\end{equation*}
In particular, the sequences $\{T_k\}_{0\leq k\leq n}$ and
$\{\mathcal{X}_k'\}_{0\leq k\leq n}$ are independent.
\item The sequence $\{\mathcal{X}_k'\}_{0\leq k\leq n}$ is a Markov chain with
  transition probabilities
\begin{equation*}
  \pP\left(\mathcal{X}_{l+1}'= {\bf s}^{i\oplus j\oplus k}\ |\ \mathcal{X}_l'={\bf s}\right) = \frac{s_i+s_j+s_k+3}{\alpha(l+1)},
\end{equation*}
where $0\leq l < n$, $1\leq i<j<k\leq N-2l$, and ${\bf
  s}=(s_1,\ldots,s_{N-2l})\in\osec$ is a generic finite
configuration with total mass $s_1+\ldots+s_{N-2l}=M$ such that
$\pP(\mathcal{X}_l'={\bf s}) > 0$.
\end{enumerate}
\end{proposition}
\begin{prooof}
  Let $0\leq l < n$, and put $L = N -2l$. By construction, the time
  $\Delta_{l+1}$ between the $l$th and the $(l+1)$th coagulation given
  $\mathcal{X}'_{l} = {\bf s} = (s_1,\ldots,s_{L})$ is exponentially
  distributed with parameter
\begin{eqnarray*}
\lefteqn{\sum_{1\leq i<j<k\leq L}(s_i +s_j +s_k +3)}\\ 
&=& 3{L\choose 3}
+\frac{1}{6}\left(\sum_{i,j,k=1}^{L}(s_i+s_j+s_k) - 3\sum_{i=1}^{L}s_i
  -3\sum_{i,j=1,\atop i\neq j}^{L}(2s_i+s_j)\right)\\
&=&\frac{1}{2}(M+L)(L-1)(L-2) = \alpha(l+1).
\end{eqnarray*}
Therefore, the waiting times
$\{\Delta_k\}_{1\leq k\leq n}$ do not depend on the states
$\{\mathcal{X}'_k\}_{1\leq k\leq n}$. The rest follows from the construction
of our process.
\end{prooof}
We turn to a description of the semigroup. Recall
that $\mathcal{X}$ starts from $\mathcal{X}(0) ={\bf r} =
(r_1,\ldots,r_N)$. In the following, $\Gamma$ denotes the Gamma function.
\begin{proposition}
\label{tc-semigroup}
  In the notation above, consider a partition $\pi$ of $\{1,\ldots,N\}$
  into $N-2l$ (non-empty) blocks $B_1,\ldots,B_{N-2l}$, each of odd
  cardinality. Denote by $\Lambda'_{\pi}(N-2l)$ the event that the $N-2l$
  atoms of $\mathcal{X}'_l$ result from the coagulation of particles $\{r_i:i\in
  B_j\}$, $j=1,\ldots,N-2l$. Then, with ${\bf r}_{B_j} = \sum_{i\in B_j}r_i$,
\begin{equation*}
  \pP\left(\Lambda'_\pi(N-2l)\right) = 
\frac{l!}{\alpha(1)\cdots\alpha(l)}\prod_{j=1}^{N-2l}\frac{\Gamma\left(({\bf
      r}_{B_j}+|B_j|+2)/2\right)(|B_j|-1)!}{\Gamma\left(({\bf r}_{B_j}+3)/2\right)((|B_j|-1)/2)!}.
\end{equation*}
\end{proposition}
\begin{prooof}
The first coagulation involves three particles with labels in the block
$B_j$ with probability
\begin{equation*}
\sum_{i<i'<i''\in B_j} \frac{r_i+r_{i'}+r_{i''}+3}{\alpha(1)} =
\frac{({\bf r}_{B_j}+|B_j|)(|B_j|-1)(|B_j|-2)}{2\alpha(1)}.
\end{equation*}
Now consider an arbitrary sequence $(k_1,\ldots,k_l)$ taking values in
$\{1,\ldots,N-2l\}$ such that for every $j=1,\ldots,N-2l$, $|\{i\leq l :
k_i = j\}|= (|B_j|-1)/2$. Using the Markov property of $\mathcal{X}'$, we see that
the probability that for all $i=1,\ldots,l$, the $i$th coagulation
affected only particles formed from initial particles with labels in
$B_{k_i}$ equals
\begin{equation*}
  \frac{1}{\alpha(1)\cdots\alpha(l)}
    \prod_{j=1}^{N-2l}\frac{\Gamma\left(({\bf r}_{B_j}+|B_j|+2)/2\right)}{\Gamma\left(({\bf
          r}_{B_j}+3)/2\right)}(|B_j|-1)!\,.
\end{equation*}
Observe that the number of such sequences $(k_1,\ldots,k_l)$ is
\begin{equation*}
{l \choose (|B_1|-1)/2,\ldots,(|B_{N-2l}|-1)/2} =\frac{l!}{((|B_1|-1)/2)!\cdots((|B_{N-2l}|-1)/2)!}.
\end{equation*}
This proves the statement.
\end{prooof}
In the setting of the proposition, denote by $\Lambda_\pi(t)$ the event
that $\mathcal{X}(t)$ has $N-2l$ atoms, each resulting from the merging of
$\{r_i : i\in B_j\}$, $j=1,\ldots,N-2l$. Since the sequence of coagulation times and the skeleton
chain $\mathcal{X}'$ are independent,
\begin{equation*}
  \pP\left(\Lambda_\pi(t)\right) =  \pP\left(T_l\leq t <
    T_{l+1},\,\Lambda'_\pi(N-2l)\right)=\pP\left(\#(t) =
    N-2l\right)\pP\left(\Lambda'_\pi(N-2l)\right). 
\end{equation*}
In particular, the semigroup of $\mathcal{X}$ is described by the
preceding proposition and the distribution of the number of
particles at time $t$, which is computed in the following lemma.
\begin{lemma}
\label{tc-nrofparticles}
In the notation above, for $l=0,\ldots,n$ and $t\geq 0$,
\begin{equation*}
  \pP\left(\#(t) = N-2l\right) =
  \sum_{j=1}^{l+1}\frac{\alpha(j)e^{-\alpha(j)t}}{\alpha(l+1)}\prod_{k=1,k\neq
      j}^{l+1}\frac{\alpha(k)}{\alpha(k)-\alpha(j)}. 
\end{equation*}
\end{lemma}
\begin{prooof}
We use
\begin{equation*}
\pP\left(\#(t) = N-2l\right) = \pP\left(T_{l+1} > t\right) -\pP\left(T_{l}>t\right).
\end{equation*}
Note that $T_k$ is distributed according to
$\sum_{i=1}^k{\alpha(i)}^{-1}{\bf e}_i$, where $\alpha(i)$ is as in the
statement of Proposition~\ref{tc-generalprop}, and ${\bf e}_1,{\bf
  e}_2,\ldots$ is a sequence of independent standard exponential variables.
As a general fact, a sum of $k$ independent exponential variables with
pairwise distinct parameters $\alpha(i)>0$ follows the hypoexponential
distribution, that is the probability distribution with density
\begin{equation*}
f(x) = \sum_{i=1}^k\alpha(i)e^{-\alpha(i)x}\prod_{j=1,j\neq i}^k\frac{\alpha(j)}{\alpha(j)-\alpha(i)}.
\end{equation*}
Integrating the density and regrouping terms result in the statement of the
lemma.
\end{prooof}

\subsection{The monodisperse case}
\label{tc-mono}
We turn to the situation where $\mathcal{X}(0)={\bf
  r}=(1,\ldots,1)$, that is the coalescent process is started from
the monodisperse configuration consisting of $N=2n+1$ atoms of unit
mass. In this case, the total mass $M$ equals $N$, so the rates $\alpha(i)$
simplify to
\begin{equation}
\label{eq:tc-monorates}
\alpha(i) = (N+1-i)(N+1-2i)(N-2i).
\end{equation}  
If ${\bf s}=(s_1,\ldots,s_{m})\in\osec$ we denote by
$\gamma({\bf s})$ the number of different $m$-tuples that can be built from
the elements $s_i$ (recall that by our convention $s_i>0$). To put it into
a formula, if $\{s_{l_i}\}_{1\leq i\leq p}$ is a maximal family of pairwise disjoint non-zero elements from
the sequence ${\bf s}$, and $k_i = \left|\{j=1,\ldots,m:
s_j=s_{l_i}\}\right|$, we define
\begin{equation*}
\gamma({\bf s}) = {m\choose k_1,\ldots,k_p}.
\end{equation*}
In other words, the ranking map
\begin{equation*}
rk: \bigcup_{m=1}^\infty\mathbb{N}^m \longrightarrow \osec
\end{equation*}
which orders $(r_1,\ldots,r_m)\in\mathbb{N}^m$ decreasingly satisfies
$|rk^{-1}(\bf{s})| =\gamma(\bf{s})$ for each ${\bf s} \neq (0,\ldots)\in\osec$.
As a corollary of Proposition~\ref{tc-semigroup}, the one-dimensional
statistics for $\mathcal{X}'$ look as follows.
\begin{corollary}
\label{1dstat-cor}
Let $0\leq l\leq n$ and ${\bf s} =(s_1,\ldots,s_{N-2l})\in\osec$ with
$s_i\in\mathbb{N}$ odd for all $i$, and $s_1+\ldots +s_{N-2l}=N$. Then,
in the situation described above,
\begin{equation*}
\pP\left(\mathcal{X}'_l = {\bf s}\right) = \gamma({\bf s})
\frac{N}{N-2l}{{N\choose l}}^{-1}\, \prod_{i=1}^{N-2l}\frac{1}{s_i}{s_i\choose \frac{s_i+1}{2}}.
\end{equation*}
\end{corollary}
\begin{prooof}
  The starting configuration is given by $(r_1,\ldots,r_N)$
  with $r_i=1$ for each $i$. Thus, if $\mathcal{X}'_l$ has $N-2l$ atoms of
  the sizes $s_1\geq\ldots\geq s_{N-2l}$, then there is a partition $\pi$
  of $\{1,\ldots,N\}$ into $N-2l$ blocks $B_1,\ldots,B_{N-2l}$ of
  cardinality $|B_j| = s_j$, such that the atoms of $\mathcal{X}'_l$
  evolved from merging the particles $\{r_i : i \in B_j\}$. Denote this
  event by $\Lambda_\pi'(N-2l)$. Since
\begin{equation*}
  \alpha(1)\cdots\alpha(l) = \frac{N!(N-1)!}{(N-l)!(N-2l-1)!}\,,
\end{equation*}
we obtain from Proposition~\ref{tc-semigroup} (note that here ${\bf r}_{B_j}=|B_j|=s_j$)
\begin{equation*}
\pP\left(\Lambda_\pi'(N-2l)\right) = \frac{(N-1-2l)!}{(N-1)!}{{N\choose
    l}}^{-1}\prod_{i=1}^{N-2l}(s_i-1)!{s_i\choose \frac{s_i+1}{2}}.
\end{equation*}
The number of such partitions $\pi$ is given by
\begin{equation*}
\frac{\gamma({\bf s})}{(N-2l)!}{N\choose s_1,\ldots,s_{N-2l}}. 
\end{equation*}
By multiplying the last two expressions together, we arrive at the stated expression.
\end{prooof}
As the reader may already check at this stage,
$\mathcal{X}'_l$ has the same distribution as the decreasingly ranked
sequence of $N-2l$ independent copies $\xi_i$ of the first hitting time of $-1$
of a simple random walk, conditioned on
$\xi_1+\ldots+\xi_{N-2l}=N$ (see Section~\ref{mconmp} for a definition of
these quantities). Indeed, if $\xi_{(k)}$ denotes the $k$th order
statistic of $\xi_1,\ldots,\xi_{N-2l}$, then for ${\bf
  s}=(s_1,\ldots,s_{N-2l})\in\osec$
\begin{eqnarray*}
\lefteqn{\pP\left((\xi_{(N-2l)},\ldots,\xi_{(1)})=
  (s_1,\ldots,s_{N-2l})\,|\,\xi_1+\ldots+\xi_{N-2l}=N\right)}\\ &=&\gamma({\bf
  s})\,
\pP\left((\xi_1,\ldots,\xi_{N-2l})=(s_1,\ldots,s_{N-2l})\,|\,\xi_1+\ldots+\xi_{N-2l}=N\right),
\end{eqnarray*}
and an application of Lemma~\ref{hittingtimes-lemma} affirms that the
last expression coincides with that obtained in the corollary. The connection
between random walks and the ternary coalescent will become much
clearer in the next section.

 \section{Duality with fragmentation via random walks}
\label{section-coalfragproc}
Our intention of this section is to prove duality of the ternary
coalescent with a fragmentation process. Let us begin with an informal
description of such processes.

Conversely to the phenomenon of coagulation of particles, one often observes
in nature or science processes of fragmentation. In these systems,
particles are broken into smaller pieces as time passes. As an example, one may think
of DNA fragmentation in biology or fractures in geophysics. Just as for
coalescent processes, one needs to impose constraints on such systems to
make them mathematically tractable. First, one assumes that the process has
no memory in the sense that the future does only depend on the present
state and not on the past. Second, one supposes that a particle is entirely
characterized by its size, that is by a real number, and third, one
requires the system to fulfill the branching property, which means that
particles split independently. 

Naively, one might first guess that a coalescent process can always be
turned into a fragmentation process by reversing time. However, even though
the memoryless property is preserved under time reversal, the
branching property is typically not fulfilled. In fact, there are only few
examples known where a duality relation holds (see~\cite{BER:3} Section 7
for an overview).

In view of our informal characterization, it is natural to call a Markov
process with values in $\osec$ a {\it ternary fragmentation process}, if
each particle splits at a certain rate according to some dislocation law
into three smaller pieces, where both the rate and the dislocation law
depend only on the particle size $s$, and the sizes of the newly formed
elements sum up to $s$. Ranked in decreasing order, these three particles
together with the ones that did not split form the next state of the
process. In particular, different particles split independently.

For our ternary coalescent starting from $N=2n+1$ atoms of
unit mass, we shall prove
\begin{theorem}
\label{thm-dual}
Reversing the coalescent chain $\{\mathcal{X}_k'\}_{0\leq k\leq n}$ in
time results in the state chain of the fragmentation process, whose
dynamics are given in Proposition~\ref{frag-proc-prop}.
\end{theorem}
We will derive our result from an explicit construction of the skeleton
chain $\mathcal{X}'$ in terms of (lengths of) excursion intervals of a
conditioned random walk. This representation will also be useful
for studying asymptotic properties in the last section.

\subsection{From configurations to paths to mass partitions}
\label{mc-and-fcts}
We first show how subsets of $\{0,1,\ldots,2n\}$ can be identified
with certain paths of nearest neighbor walks on $\mathbb{Z}$ of length
$2n+1$. The excursion intervals above two consecutive (new) minima of such
paths partition the space $\mathbb{Z}/(2n+1)\mathbb{Z}$ into discrete
arcs. Taking the ranked sequence of their lengths, we obtain the main
object of our interest.

To begin with, define the configuration space $\cs_n$ to be the set of all
subsets of $\{0,\ldots,2n\}$ which have cardinality less or equal to $n$. We
often represent $x\in\cs_n$ by the vector $(x(i))_{0\leq i\leq 2n}$, where
\begin{equation*}
x(i) = \left\{\begin{array}{l@{\ ,\ \ }l}
      1 & i\in x\\
      0 & i\notin x\end{array}\right..
\end{equation*}
Under this identification, we may regard $x$ as a mass distribution. We use
the terminology that a site $i$ is occupied by a mass if $x(i) = 1$ and vacant otherwise.
The number of occupied sites (the cardinality of the subset $x$) is denoted by
\begin{equation*}
|x| =  \left|\{i\in\{0,\ldots,2n\}: x(i) =1\}\right|.
\end{equation*}
We identify a configuration $x\in \cs_n$ with a path of a nearest neighbor walk of
length $2n+1$ on $\mathbb{Z}$ in the following way. Starting from the origin at time zero,
the walk goes one step up if site $0$ is occupied, i.e. $x(0)=1$, and
down otherwise, then above if $x(1)=1$, down if $x(1) = 0$ and so on, up to
time $2n$. More precisely, the corresponding path $S(x)$ is given by
$S(x)_0 = 0$ and for $1\leq j\leq 2n+1$,
\begin{equation*}
S(x)_j =  
2\left(\sum_{i=0}^{j-1}x(i)\right) - j.
\end{equation*}
Notice that by definition, $S(x)_{2n+1}= 2(|x|-n)-1$.  
Clearly, the mapping $\cs_n\ni x\mapsto S(x)$ is one-to-one.

As we show next, the excursion intervals of such a path provide us
with an element $\vp_1(x)$ in the space of cyclically ordered partitions of
$\mathbb{Z}/(2n+1)\mathbb{Z}$ into discrete arcs,
\begin{eqnarray*}
  \omp_{2n+1} &=& \left\{{\tt s^\circ} =({\tt s}_1,\ldots,{\tt s}_m) :
    \mbox{there exist }a_1<a_2<\ldots<a_m\leq 2n+1,\right.\\
  &&\ \: m,a_i\in\mathbb{N},\mbox{ such that for }1\leq i\leq m-1, {\tt s}_i =
  [a_i,a_{i+1})\cap\mathbb{N},\\ 
  &&\,\:  \left.{\tt s}_m = \left([a_m, 2n+1)\cup[0,a_1)\right)\cap\mathbb{Z}_+\right\}.
\end{eqnarray*}
Take $x\in\cs_n$, and let $M= -S(x)_{2n+1}$. With $\underline{m}(x) = \min_{0\leq j\leq
  2n+1}S_j(x)$, define the first time at which $S(x)$
reaches $\underline{m}(x) +k$, $k=0,\ldots,M-1$,
\begin{equation*}
m_k(S(x)) = \inf \left\{j\geq 0 : S_j(x) = \underline{m}(x) +k\right\}.
\end{equation*}
For $i=1,\ldots,M$, put $a_i = m_{M-i}(S(x))$.  We 
construct a sequence ${\tt s^\circ}=\left({\tt s}_1,\ldots,{\tt s}_M
\right)\in\omp_{2n+1}$ by setting ${\tt s}_i=[a_i,a_{i+1})\cap\mathbb{N}$
for $i=1,\ldots,M-1$, ${\tt s}_M = \left([a_M,
  2n+1)\cup[0,a_1)\right)\cap\mathbb{Z}_+$. In other words, if we look for $k=0,\ldots,2n$ at the shifted path
$\theta_k(S(x))$ defined by
\begin{equation*}
\theta_k(S(x))_i = \left\{\begin{array}{l@{\ ,\ \ }l}
     S(x)_{i+k}-S(x)_k  & 0\leq i\leq 2n+1-k\\
     S(x)_{i+k-(2n+1)} +S(x)_{2n+1}-S(x)_k &
     2n+1-k<i\leq 2n+1\end{array}\right.,
\end{equation*}
then the element ${\tt s^\circ}$ corresponds to the $M$ successive excursion intervals of
$\theta_{m_{M-1}}S(x)$ above two consecutive (new) minima. 
The length $|{\tt s}_i|$ of such an interval is also referred to as
a {\it ladder epoch}. We let $\vp_1(x) = {\tt
  s^\circ}$ and define $\vp_2$ as the function which sends ${\tt s^\circ} =({\tt
  s}_1,\ldots,{\tt s}_m)\in\omp_{2n+1}$ to its arc lengths
$\{|{\tt s}_i|\}_{1\leq i\leq m}$, arranged in decreasing order. In this way, we
obtain an element in the space of mass partitions
\begin{equation*}
  \ump_{2n+1}= \left\{{\bf s} =(s_1,\ldots,s_m) : s_1\geq
    s_2\geq\ldots\geq s_m\,, m,s_i \in\mathbb{N},\, \sum_{i=1}^ms_i=
    2n+1\right\}.
\end{equation*}
By filling up with an infinite sequence of zeros, we will often identify mass
partitions with elements in $\osec$.  To summarize our construction, the
concatenation map $\vp$
\begin{equation*}
 \vp = \vp_2\circ\vp_1 : \cs_n\stackrel{\vp_1}\longrightarrow \omp_{2n+1}\stackrel{\vp_2}\longrightarrow\ump_{2n+1}\subset\osec.
\end{equation*}
sends configurations $x\in\cs_n$ via their path representations to
partitions of $\mathbb{Z}/(2n+1)\mathbb{Z}$ and then to mass
partitions.
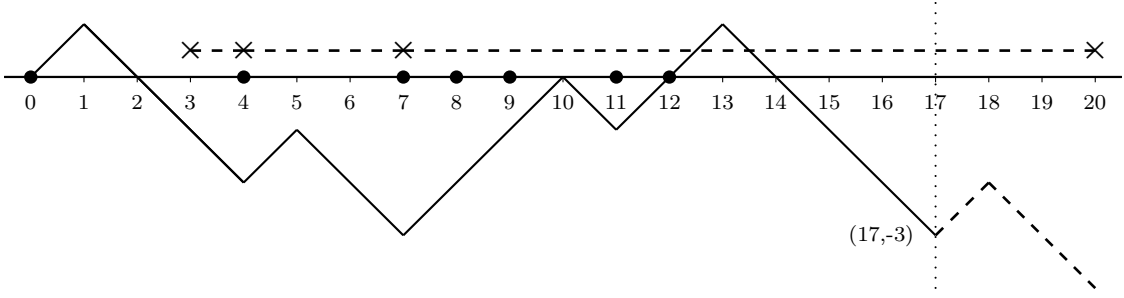
\begin{figure}
\psset{xunit=0.7cm,yunit=0.7cm,algebraic=true,dotstyle=o,dotsize=3pt
      0,linewidth=0.8pt,arrowsize=3pt 2,arrowinset=0.25}
    \centering
    \begin{pspicture*}(-0.5,-4.5)(20.5,1.5) 
      \psaxes[labelFontSize=\scriptstyle,xAxis=true,yAxis=false,labels=x,Dx=1,Dy=1,ticksize=-2pt
      0,subticks=0]{-}(0,0)(-0.5,-4.5)(20.5,1.5) \psline(0,0)(1,1)
      \psline(1,1)(4,-2) \psline(1,1)(4,-2) \psline(5,-1)(4,-2)
      \psline(5,-1)(7,-3) \psline(7,-3)(10,0) \psline(10,0)(11,-1)
      \psline(13,1)(11,-1) \psline(13,1)(17,-3)
      \psline[linewidth=1.0pt,linestyle=dashed,dash=4pt 4pt](3,0.5)(4,0.5)
      \psline[linewidth=1.0pt,linestyle=dashed,dash=4pt 4pt](4,0.5)(7,0.5)
      \psline[linewidth=1.0pt,linestyle=dashed,dash=4pt 4pt](17,-3)(18,-2)
      \psline[linewidth=1.0pt,linestyle=dashed,dash=4pt 4pt](18,-2)(20,-4)
      \psline[linestyle=dotted](17,2)(17,-4)
      \psline[linewidth=1.0pt,linestyle=dashed,dash=4pt 4pt](7,0.5)(20,0.5)
 \begin{scriptsize}
  \rput[Tr](16.6,-3.0){(17,-3)}
   \psdots[dotsize=5pt 0,dotstyle=*](0,0) \psdots[dotsize=5pt
   0,dotstyle=*](4,0) \psdots[dotsize=5pt 0,dotstyle=*](7,0)
   \psdots[dotsize=5pt 0,dotstyle=*](8,0) \psdots[dotsize=5pt
   0,dotstyle=*](9,0) \psdots[dotsize=5pt 0,dotstyle=*](11,0)
   \psdots[dotsize=5pt 0,dotstyle=*](12,0) \psdots[dotsize=7pt
   0,dotstyle=x](3,0.5) \psdots[dotsize=7pt 0,dotstyle=x](4,0.5)
   \psdots[dotsize=7pt 0,dotstyle=x](7,0.5) \psdots[dotsize=7pt
   0,dotstyle=x](20,0.5)
 \end{scriptsize}
\end{pspicture*}
\caption{The black dots represent the configuration
  $x=\{0,4,7,8,9,11,12\}\subset \cs_8$. The corresponding path $S(x)$
  starts at zero and ends in $-3$ at time $17$. It is periodically extended
  up to time $20$ to better recognize the excursion intervals
  $\vp_1(x)$. They are visualized by the dashed line above the $x$-axis,
  where the crosses mark the endpoints of the intervals,
  i.e. $\vp_1(x) = ({\tt s}_1,{\tt s}_2,{\tt s}_3)$ with ${\tt s}_1=[3,4)\cap\mathbb{N}$,
  ${\tt s}_2=[4,7)\cap\mathbb{N}$, ${\tt s}_3=([7,17)\cup[0,3))\cap\mathbb{Z}_+$.}
\end{figure}

\subsection{Random evolution}
\label{randomevolution}
Our purpose here is to randomize the input of the map $\vp:
\cs_n\rightarrow\ump_{2n+1}$ to obtain (a sequence of) random mass
partitions. More precisely, we construct two Markov chains on
$\cs_n$ running from time zero up to $n$ as follows. Let $X=\{X_k\}_{0\leq k\leq n}$ be the Markov chain with $X_0 = \emptyset$
and transition probabilities
\begin{equation*}
p_X(x,y) = \left\{\begin{array}{c@{\ ,\ \ }l}
      \frac{1}{2n+1-|x|} & x\subset y\mbox{ and } y\backslash x = \{i\}\mbox{
        for some }i\in\{0,\ldots,2n\}\backslash x\\
      0 & \mbox{otherwise}\end{array}\right..
\end{equation*}
In words, $(X_0,\ldots,X_l)$ is obtained by occupying successively $l$
sites from $\{0,\ldots,2n\}$, chosen uniformly at random. From the point of
view of sets, $X_l$ is uniformly distributed on the space of all
$l$-subsets of $\{0,\ldots,2n\}$. By identifying with the random path
$S(X_l)$, we will also think of $X_l$ as simple random walk up to
time $2n+1$, conditioned to end at position $-2(n-l)-1$.

Let $Y=\{Y_k\}_{0\leq k\leq n}$ be the Markov chain with $Y_0$ being uniformly
distributed on the space of all $n$-subsets of $\{0,\ldots,2n\}$ and
transition probabilities
\begin{equation*}
p_Y(x,y) = \left\{\begin{array}{c@{\ ,\ \ }l}
      \frac{1}{|x|} & y\subset x\mbox{ and } x\backslash y = \{i\}\mbox{
        for some }i\in\{0,\ldots,2n\}\backslash y\\
      0 & \mbox{otherwise}\end{array}\right..
\end{equation*}
In words, $(Y_0,\ldots,Y_l)$ is obtained by removing successively $l$
masses chosen uniformly at random from the starting configuration $Y_0$.
In terms of sets, $Y_l$ is uniformly distributed on the space of all
$(n-l)$- subsets of $\{0,\ldots,2n\}$.  Similarly as above, $Y_l$ can be
identified with simple random walk up to time $2n+1$, conditioned to end at
$-2l-1$. Note that by construction, we have the duality relation
\begin{equation}
\label{eq:dualityXY}
(X_0,\ldots,X_n) \eql (Y_n,\ldots,Y_0).
\end{equation}
\subsection{Realization of the skeleton chains}
\label{mconmp}
We are not interested in $X$ and $Y$ themselves, but rather in $\vp(X) =
\{\vp(X_k)\}_{0\leq k\leq n}$ and $\vp(Y)$. As we will show in
Proposition~\ref{coal-proc-prop}, the former is the state chain of the
ternary coalescent starting from $N=2n+1$ atoms of unit mass. The latter
is characterized by Proposition~\ref{frag-proc-prop} as the state chain
of a fragmentation process starting from a single particle of mass $N$.

We need some preparation. Recall that simple random walk on
$\mathbb{Z}$ is the Markov chain $S = \{S_m\}_{m\geq 0}$ with $S_0 = 0$ and
$S_m = \zeta_1+\ldots+\zeta_m$, where $\zeta_1,\zeta_2,\ldots$ are
independent random variables with $\pP(\zeta_i =\pm 1)= 1/2$. For
$k\in\mathbb{Z}$, the first hitting time of $k$ is
denoted by
\begin{equation*}
H_k= \inf \{m\geq 1: S_m= k\}.
\end{equation*}
The following result on the distribution of $H_k$ is classical.
\begin{lemma}
\label{hittingtimes-lemma}
Let $k\in\mathbb{Z}$, $k\neq 0$, and $m\in\mathbb{N}$. Then
\begin{equation*}
  \pP\left(H_k = m\right) = \left\{\begin{array}{l@{\ ,\ \ }l}
      \frac{|k|}{m}{m\choose(m+|k|)/2}2^{-m} & k=m[\textup{mod }2]\\
      0 & k\neq m[\textup{mod }2]\end{array}\right..
\end{equation*}
Moreover, if $m=2n+1$ and $k$ is a fixed odd number, as $n\rightarrow \infty$, 
\begin{equation*}
  \pP\left(H_k = m\right) \sim \frac{1}{2}\sqrt{\frac{1}{\pi n^3}}.
\end{equation*}
\end{lemma}
\begin{prooof}
  Clearly, for the probability to be different from zero the numbers $k$ and $m$ must
  have the same parity. Then, using the hitting time theorem (see for
  example~\cite{KMP}) in the first equality,
\begin{equation*}
\pP(H_k = m) = \frac{|k|}{m}\pP(S_m=k)=\frac{|k|}{m}{m\choose(m+|k|)/2}2^{-m}. 
\end{equation*}
The second statement follows from Stirling's formula for the factorial.
\end{prooof}
Before looking at $\vp(X)$ and $\vp(Y)$ in detail, let us give an
indication that the former is the skeleton chain of the ternary
coalescent. Recall Corollary~\ref{1dstat-cor} and the connection between
$\mathcal{X}'$ and hitting times. Let $N=2n+1$, $0\leq l\leq n$,
and take $N$ independent copies $\xi_i$ of the hitting time
$H_{-1}$. Denote by $\xi_{(k)}$ the $k$th order statistic of
$\xi_1,\ldots,\xi_{N-2l}$.
\begin{proposition}
\label{1dstat-prop-RW}
$\vp(X_l)$ is distributed according to
$(\xi_{(N-2l)},\ldots,\xi_{(1)})$ conditionally on
$\xi_1+\ldots+\xi_{N-2l}=N$, i.e. the one-dimensional
distributions of $\vp(X)$ and
$\mathcal{X}'$ started from $N$ atoms of mass one agree. 
\end{proposition}
\begin{prooof}
  We identify $X_l$ with simple random walk $S(X_l)$ up to time $N$,
  conditioned to end at $-(N-2l)$.  For notational simplicity, let us write
  $S$ instead of $S(X_l)$. Also recall the definitions of 
  $\theta_k(S)$ and $m_k(S)$ from Section~\ref{mc-and-fcts}.  By Theorem 1
  of~\cite{BCP}, if $\nu$ is a uniform random variable on
  $\{0,\ldots,N-2l-1\}$ independent of $S$, then the chain
  $\theta_{m_\nu}(S)$ has the law of $S$ conditioned on
  $H_{-(N-2l)}=N$. Moreover, the index $m_{\nu}$ is uniformly distributed
  on $\{0,\ldots,N-1\}$ and independent of the chain
  $\theta_{m_\nu}(S)$. Denote by $\theta_kX_l$ the shifted configuration
  defined by $\theta_kX_l(i) = X_l(i+k[\textup{mod }N])$. Clearly,
  $\vp(X_l) = \vp(\theta_kX_l)$ for each $k$. From Theorem 1 of~\cite{BCP}
  we thus infer that for $(s_1,\ldots,s_{N-2l})\in\osec$,
\begin{eqnarray*}
\lefteqn{\pP\left(\vp(X_l) =
  (s_1,\ldots,s_{N-2l})\right)=\pP\left(\vp(\theta_{m_\nu}X_l) =
(s_1,\ldots,s_{N-2l})\right)}\\ 
&=& \pP\left((\xi_{(N-2l)},\ldots,\xi_{(1)})=
(s_1,\ldots,s_{N-2l})\,|\,\xi_1+\ldots+\xi_{N-2l}=N\right).
\end{eqnarray*}
\end{prooof}
For the moment, we leave $\vp(X)$ aside and first turn to $\vp(Y)$. In the
sequel it is convenient to use the notion of multisets, which we
distinguish from normal sets by using double braces. For example,
$\{\{a,b,c,c\}\}$ contains the elements $a$, $b$ each with multiplicity $1$
and the element $c$ with multiplicity $2$. The cardinality of this multiset
is $4$, the order of elements is irrelevant, as for sets.

Let $\xi_1,\xi_2,\xi_3$ be three
independent copies of the hitting time $H_{-1}$. To state the transition
mechanism of $\vp(Y)$ in a concise way, we define a family $\mu=
\left(\mu_s,s\geq 3\mbox{ odd}\right)$ of probability laws, supported on
\begin{equation*}
\Omega_s = \left\{R=\{\{r_1,r_2,r_3\}\} : r_i\in\mathbb{N}\mbox{ odd},\, 
r_1+r_2+r_3 = s\right\},
\end{equation*}
by setting
\begin{equation}
\mu_s(R) = \pP\left(\{\{\xi_1,\xi_2,\xi_3\}\}=R\ |\ \xi_1+\xi_2+\xi_3 = s\right).
\end{equation}
More explicitly, applying Lemma~\ref{hittingtimes-lemma} results in the
expression
\begin{equation}
\label{eq:disloclaw}
\mu_s(R) =
\gamma\frac{s}{3r_1r_2r_3}{r_1\choose\frac{r_1+1}{2}}{r_2\choose\frac{r_2+1}{2}}{r_3\choose\frac{r_3+1}{2}}\left[{s\choose\frac{s+3}{2}}\right]^{-1},
\end{equation}
where $\gamma$ is the number of triplets $(r_i,r_j,r_k)$ that can be formed
from $R=\{\{r_1,r_2,r_3\}\}$,
\begin{equation*}
\gamma =\left\{\begin{array}{l@{\ ,\ \ }l}
       6 & |\{r_1,r_2,r_3\}| = 3\\
       3 & |\{r_1,r_2,r_3\}| = 2\\
       1 & |\{r_1,r_2,r_3\}| = 1\end{array}\right..
\end{equation*}
\begin{proposition}
\label{frag-proc-prop}
  $\vp(Y) = \{\vp(Y_k)\}_{0\leq k\leq n}$ is a Markov chain. Its transition
  mechanism from time $l\leq n-1$ to $l+1$ is described as follows. 
  \begin{enumerate}
  \item[(a)] Conditionally on $\vp(Y_l) =
    {\bf s}=(s_1,\ldots,s_{2l+1})\in\ump_{2n+1}$, select an index
    $\iota\in\{1,\ldots,2l+1\}$ according to the law
    \begin{equation*}
  \pP\left(\iota = i\ |\ \vp(Y_l) = {\bf s}\right) = \frac{s_i-1}{2(n-l)}.
    \end{equation*}
  \item[(b)] Given $\vp(Y_l) =
    {\bf s}$ and $\iota=i$, split $s_i$ according to the law $\mu_{s_i}$ into
    three numbers and rank them together with $s_m$, $m\in\{1,\ldots,2l+1\}\backslash\{i\}$,
    in decreasing order to obtain a new mass partition.
  \end{enumerate}
\end{proposition}
\begin{prooof}
  Fix $l\in\{0,\ldots,n-1\}$. We write $\vp(Y)_{0:i}$ for the vector
$\left(\vp(Y_0),\ldots,\vp(Y_i)\right)$.
The Markov property will follow from
  \begin{enumerate}
  \item $\vp(Y)_{0:l}$ and $\vp(Y_{l+1})$ are conditionally independent
    given $\vp_1(Y_l)$.
    \item \label{eq:mpY-1}
    $\vp_1(Y_l)$ and $\vp(Y_{l+1})$ are conditionally independent given $\vp(Y_l)$.
  \end{enumerate}
 Indeed, assuming (i) and (ii), we have for
  ${\bf r}_{0:l+1}=({\bf r}_0,\ldots,{\bf r}_{l+1})\in\ump_{2n+1}\times\ldots\times\ump_{2n+1}$,
\begin{eqnarray*}
  \lefteqn{\pP\left(\vp(Y)_{0:l+1} = {\bf r}_{0:l+1}\right)}\\
  &=&\sum_{{\tt u}:
    \vp_2({\tt u^\circ})={\bf r}_l}\pP\left(\vp(Y)_{0:l}={\bf r}_{0:l}\ |\ \vp_1(Y_l) =
    {\tt u^\circ}\right)\pP\left(\vp(Y_{l+1})={\bf r}_{l+1}\ |\ \vp_1(Y_l) =
    {\tt u^\circ}\right)\\
&&\times\ \pP\left(\vp_1(Y_l)={\tt u^\circ}\right)\\
  & =& \pP\left(\vp(Y_{l+1}) = {\bf r}_{l+1}\ |\ \vp(Y_l) =
  {\bf r}_l\right)\pP\left(\vp(Y)_{0:l}={\bf r}_{0:l}\right). 
\end{eqnarray*}
For (i), the key step is to show that the conditional law of $\vp(Y)_{0:l}$ given $Y_l$ only
depends on $\vp_1(Y_l)$. In that direction, we work conditionally on
$\vp_1(Y_l)={\tt s^\circ} = ({\tt s}_1,\ldots,{\tt s}_{2l+1})$ and denote
by $N_k(i)= |Y_k\cap {\tt s}_i|$ the number of sites of the arc ${\tt s}_i$ which are occupied by $Y_k$. Write
$N_k$ for the family $\{N_k(i)\}_{1\leq i\leq 2l+1}$. Let $i_l$ denote the
unique index such that the singleton $Y_{l-1}\backslash Y_l \subset {\tt
  s}_{i_l}$. In other words, $i_l$ is the unique index $i$ such that
$N_{l-1}(i) =N_l(i)+1$. Then $\vp_1(Y_{l-1})$ results from $\vp_1(Y_l)={\tt
  s}$ by merging the arcs ${\tt s}_{i_l}$, ${\tt s}_{i_l+1}$ and ${\tt
  s}_{i_l+2}$ (with the convention that indices of arcs are taken modulo
$2l+1$).  By iteration, we realize that the sequence
$N_{0:l}=(N_0,\ldots,N_l)$ determines $\vp_1(Y)_{0:l}$ and therefore
also $\vp(Y)_{0:l}$. Hence it now suffices to check that the conditional
distribution of $N_{0:l}$ given $Y_l$ only depends on $\vp_1(Y_l) = {\tt
  s^\circ}$, which is straightforward from the dynamics and the observation that
for every $i=1,\ldots,2l+1$, the arc ${\tt s}_i$ has exactly $(|{\tt s}_i|
+1)/2$ sites which are not occupied by $Y_l$.

We are now able to prove (i). Take
${\bf t}\in\ump_{2n+1}$ with $\pP\left(\vp_1(Y_l) = {\tt
    s^\circ},\vp(Y_{l+1}) = {\bf t}\right) >0$. Then
\begin{eqnarray*}
\lefteqn{\pP\left(\vp(Y)_{0:l} = {\bf r}_{0:l}, \vp_1(Y_l) = {\tt s^\circ}, 
\vp(Y_{l+1}) = {\bf t}\right)} \\
&=& \sum_{x\in\vp_1^{-1}({\tt s^\circ}),\atop y\in\vp^{-1}({\bf t})}\pP\left(\vp(Y)_{0:l}
  = {\bf r}_{0:l}\ | \ Y_l=x, Y_{l+1} = y\right)\pP\left(Y_l=x, Y_{l+1} =
  y\right).
\end{eqnarray*} 
Since $Y$ is a Markov chain, it follows that for $x,y\in\cs_n$ with
$\pP(Y_l = x, Y_{l+1} = y)>0$,
\begin{equation*}
\pP\left(\vp(Y)_{0:l} = {\bf r}_{0:l}\ | \ Y_l=x, Y_{l+1} = y\right) = 
\pP\left(\vp(Y)_{0:l} = {\bf r}_{0:l}\ | \ Y_l=x\right).
\end{equation*}
Plugging this into the above formula and using the conditional independence
of $\vp(Y)_{0:l}$ and $Y_l$ given $\vp_1(Y_l)$, we deduce
that for $x\in\vp_1^{-1}({\tt s^\circ})$, 
\begin{equation*}
\pP\left(\vp(Y)_{0:l} = {\bf r}_{0:l}\ |\ \vp_1(Y_l) = {\tt s^\circ},
    \vp(Y_{l+1}) = {\bf t}\right) = \pP\left(\vp(Y)_{0:l} = {\bf r}_{0:l}\ | \ Y_l=x\right).
\end{equation*}
Similarly, one sees that the right hand side equals $\pP(\vp(Y)_{0:l} =
{\bf r}_{0:l}\ |\ \vp_1(Y_l) = {\tt s^\circ})$, and (i) follows. We turn to
(ii) and the description of the transition mechanism. 
We keep the conditioning on $\vp_1(Y_l)={\tt s^\circ}$. Note that $Y_{l+1}$
evolves from $Y_l$ by removing uniformly at random one of the $n-l$
masses. By identifying $Y_l$ with $S(Y_l)$,
this amounts to switching one of the upward steps chosen uniformly at
random into a downward step. More precisely, under our conditioning, the
$i$th arc ${\tt s}_i$ is picked with probability
\begin{equation}
\label{eq:probaitharc}
\frac{\mbox{number of upward steps over }{\tt s}_i}{\mbox{total number of
    upward steps}} = \frac{(|{\tt s}_i|-1)/2}{n-l}, 
\end{equation}
then one of the upward steps over ${\tt s}_i$ is selected with uniform
probability and changed into a downward step.
Up to a vertical shift in space, $S(Y_l)$ restricted to the arc ${\tt s}_i$ obeys the law of
simple random walk conditioned on $H_{-1}=|{\tt s_i}|$ (with an obvious
modification for the last arc ${\tt s}_{2l+1}$). Given an upward
step over ${\tt s_i}$ is switched, $S(Y_{l+1})$ restricted to
${\tt s}_i$ can therefore be seen as simple random walk conditioned on
$H_{-3}= |{\tt s_i}|$. In terms of $\vp(Y)$, we deduce that 
$\vp(Y_{l+1})$ is obtained by first picking the $i$th arc ${\tt s}_i$ with 
probability given in~\eqref{eq:probaitharc}, then splitting its length
according to $\mu_{|\tt{s}_i|}$ into three numbers $r_1,
r_2, r_3$ corresponding to the first three ladder epochs of simple random
walk conditioned on $H_{-3}=|{\tt s}_i|$, and finally ranking them together with the numbers
$|{\tt s}_j|$, $j\neq i$, in decreasing order. In particular, we realize that for predicting $\vp(Y_{l+1})$ out of
$\vp_1(Y_l)$, the additional information given by $\vp_1(Y_l)$ compared to
$\vp(Y_l)$, namely the location of the arcs, is irrelevant. Hence also
(ii) holds.
\end{prooof}

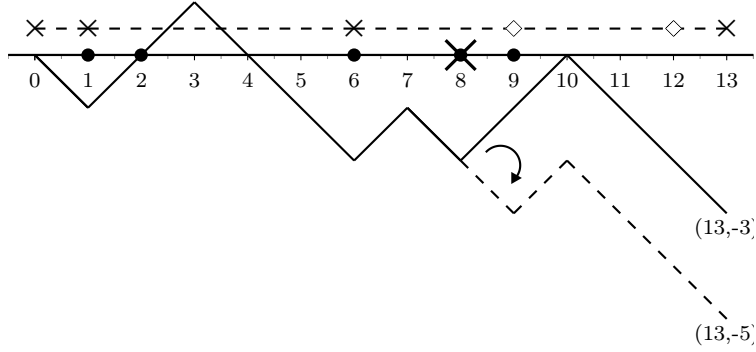
\begin{figure}
\centering
\psset{xunit=0.7cm,yunit=0.7cm,algebraic=true,dotstyle=o,dotsize=3pt 0,linewidth=0.8pt,arrowsize=3pt 2,arrowinset=0.25}
\begin{pspicture*}(-0.5,-5.5)(13.5,1.5)
  \psaxes[labelFontSize=\scriptstyle,xAxis=true,yAxis=false,Dx=1,Dy=1,ticksize=-2pt
  0,subticks=2]{-}(0,0)(-0.5,-5.5)(13.5,1.5) \psline(0,0)(1,-1)
  \psline(0,0)(1,-1) \psline(1,-1)(3,1) \psline(3,1)(6,-2)
  \psline(8,-2)(7,-1) \psline(8,-2)(7,-1) \psline(8,-2)(10,0)
  \psline(10,0)(13,-3) \psline(6,-2)(7,-1)
  \psline[linestyle=dashed,dash=4pt 4pt](0,0.5)(1,0.5)
  \psline[linestyle=dashed,dash=4pt 4pt](1,0.5)(6,0.5)
  \psline[linestyle=dashed,dash=4pt 4pt](6,0.5)(13,0.5)
  \psline[linestyle=dashed,dash=4pt 4pt](8,-2)(9,-3)
  \psline[linestyle=dashed,dash=4pt 4pt](9,-3)(10,-2)
  \psline[linestyle=dashed,dash=4pt 4pt](10,-2)(13,-5)
  \parametricplot{-0.7309067071567155}{2.4106859464330777}{1*0.37*cos(t)+0*0.37*sin(t)+8.75|0*0.37*cos(t)+1*0.37*sin(t)+-2.09}
\begin{scriptsize}
\psdots[dotsize=5pt 0,dotstyle=*](1,0)
\psdots[dotsize=5pt 0,dotstyle=*](2,0)
\psdots[dotsize=5pt 0,dotstyle=*](6,0)
\psdots[dotsize=5pt 0,dotstyle=*](8,0)
\psdots[dotsize=5pt 0,dotstyle=*](9,0)
\psdots[dotsize=7pt 0,dotstyle=x](0,0.5)
\psdots[dotsize=7pt 0,dotstyle=x](1,0.5)
\psdots[dotsize=7pt 0,dotstyle=x](6,0.5)
\psdots[dotsize=7pt 0,dotstyle=x](13,0.5)
\psdots[dotsize=5pt 0,dotstyle=square,dotangle=45](9,0.5)
\psdots[dotsize=5pt 0,dotstyle=square,dotangle=45](12,0.5)
\psdots[dotsize=13pt 0,dotstyle=x](8,0)
\psdots[dotsize=4pt 0,dotstyle=triangle*,dotangle=270](9.02,-2.33)
\rput[t](13,-3.1){(13,-3)}  
\rput[t](13,-5.1){(13,-5)}
\end{scriptsize}
\end{pspicture*}
\caption{The transition mechanism from $\vp(Y_1)$ to $\vp(Y_2)$, where
  $n= 6$. Here, at time $1$ the chain $Y$ is in the configuration state
  $Y_1=\{1,2,6,8,9\}$. Then the mass at position $8$ is removed. For the
  corresponding path, this means that the upward step at time $8$ is
  changed into a downward step. The new path $S(Y_2)$ coincides up to time
  $8$ with the old path $S(Y_1)$ and is then indicated by the dashed
  line. The excursion interval $[6,13)$ is broken into
  three intervals $[6,9)$, $[9,12)$, $[12,13)$. Therefore,
  $\vp_1(Y_2)=({\tt s}_1, {\tt s}_2,{\tt s}_3,{\tt s}_4,{\tt s}_5)$ with
  ${\tt s}_1 = [1,6)\cap\mathbb{N}$, ${\tt s}_2 = [6,9)\cap\mathbb{N}$, ${\tt s}_3 =
[9,12)\cap\mathbb{N}$, ${\tt s}_4 = [12,13)\cap\mathbb{N}$, ${\tt s}_5=[0,1)\cap\mathbb{Z}_+$ and $\vp(Y_2) = (5,3,3,1,1)$.}
\end{figure}

Let us now characterize
$\vp(X)$.
\begin{proposition}
\label{coal-proc-prop}
  $\vp(X) = \{\vp(X_k)\}_{0\leq k\leq n}$ is a Markov chain. Its transition
  mechanism from time $l\leq n-1$ to $l+1$ is described as follows. 
  \begin{enumerate}
  \item[(a)] Conditionally on $\vp(X_l) =
    {\bf s}=(s_1,\ldots,s_{2(n-l)+1})\in\ump_{2n+1}$, select an index
    $\iota$ out of the set of all $3$-subsets of $\{1,\ldots,2(n-l)+1\}$
    according to the law
    \begin{equation*}
  \pP\left(\iota = \{i,j,k\}\ |\ \vp(X_l) = {\bf s}\right) = 
     \frac{s_i+s_j+s_k+3}{(2n+1-l)2(n-l)(2(n-l)-1)}\, .
    \end{equation*}
  \item[(b)] Given $\vp(X_l)={\bf s}$ and $\iota=\{i,j,k\}$, rank the sum
    $r= s_i+s_j+s_k$ together with the numbers $s_m$,
    $m\in\{1,\ldots,2(n-l)+1\}\backslash\{i,j,k\}$, in decreasing order to
    obtain a new mass partition.
\end{enumerate}
\end{proposition}
\begin{prooof}
From the duality~\eqref{eq:dualityXY} it follows that $\vp(X)$ is obtained
by reversing $\vp(Y)$ in time. In particular, the Markov property carries
over from $\vp(Y)$ to $\vp(X)$.
 
It remains to look at the transition mechanism.
The step from $l=n-1$ to $n$ is obvious from the construction of $X$ and
$\vp$. Now fix $l\in\{0,\ldots,n-2\}$, and let $M=2(n-l)+1$. We work
conditionally on $\vp(X_l) ={\bf s}=(s_1,\ldots,s_M)\in\ump_{2n+1}$. By
construction, $\vp(X_{l+1})$ is obtained from $\vp(X_l)$ by summing up
three numbers $s_i$,$s_j$,$s_k$, where $i$,$j$,$k$ are pairwise distinct,
and rearranging the sum together with $s_m$, $m\neq i,j,k$, in decreasing
order. Write ${\tt s^\circ} = ({\tt s}_1,\ldots,{\tt s}_M)$ for the
partition $\vp_1(X_l)$, and let $\nu$ be uniformly distributed on
$\{0,\ldots,M-1\}$, independent of $X_l$. By the random walk representation
and Theorem 1 of~\cite{BCP}, the law of the cyclically ordered arc lengths
$\left(|{\tt s}_{1+\nu}|,\ldots,|{\tt s}_{M+\nu}|\right)$ (indices are
taken modulo $M$) agrees with the law of the $M$ subsequent ladder epochs
of simple random walk conditioned on $H_{-M}=2n+1$. In particular, the law
of $\left(|{\tt s}_{1+\nu}|,\ldots,|{\tt s}_{M+\nu}|\right)$ is invariant
under permutations and therefore equals the law of 
$\left(s_{\sigma(1)},\ldots,s_{\sigma(M)}\right)$, where $\sigma$ is a
permutation of $\{1,\ldots,M\}$, chosen uniformly at random and
independently of $X_l$. Note that this can also be deduced directly from
the fact that $X_l$ is uniformly distributed on the space of all
$l$-subsets of $\{0,\ldots,2n\}$.  The probability that $s_i$,$s_j$,$s_k$
are replaced by their sum is given by the probability that the arcs ${\tt
  s}_{\sigma^{-1}(i)+\nu}$, ${\tt s}_{\sigma^{-1}(j)+\nu}$, ${\tt
  s}_{\sigma^{-1}(k)+\nu}$ merge. This is the case if and only if the arcs
adjoin each other and the singleton $X_{l+1}\backslash X_l$ is contained in
that arc which is followed in clockwise order by the other two. More
formally, the arcs merge if and only if there is a permutation $\rho$ of
the indices $i,j$ and $k$ such that $X_{l+1}\backslash X_l\subset {\tt
  s}_{\sigma^{-1}(\rho(i))+\nu}$, and $\sigma^{-1}(\rho(j)) =
\sigma^{-1}(\rho(i)) +1$, $\sigma^{-1}(\rho(k)) = \sigma^{-1}(\rho(i)) +2$
(both equalities are taken modulo $M$).  Given $X_{l+1}\backslash
X_l\subset {\tt s}_{\sigma^{-1}(i)+\nu}$, the probability that ${\tt
  s}_{\sigma^{-1}(i)+\nu}$, ${\tt s}_{\sigma^{-1}(j)+\nu}$, ${\tt
  s}_{\sigma^{-1}(k)+\nu}$ merge is therefore
\begin{eqnarray*}
\frac{2}{M-1}\times\frac{1}{M-2}.
\end{eqnarray*}
The probability that
$X_{l+1}\backslash X_l\subset{\tt s}_{\sigma^{-1}(i)+\nu}$ is 
\begin{equation*}
\frac{\mbox{number of vacant sites in }{\tt s}_{\sigma^{-1}(i)+\nu}\mbox{
    at time }l}{\mbox{total number of vacant sites at time }l} =
\frac{(s_i+1)/2}{2n+1-l}.
\end{equation*}
Altogether, given $\vp(X_l) ={\bf s}$,
\begin{equation*}
  \pP\left({\tt s}_{\sigma^{-1}(i)+\nu}, {\tt s}_{\sigma^{-1}(j)+\nu}, {\tt
      s}_{\sigma^{-1}(k)+\nu} \mbox{ merge}\right) =
  \left(\frac{(s_i+s_j+s_k+3)/2}{2n+1-l}\right) 
  \frac{2}{M-1}\times\frac{1}{M-2},
\end{equation*}
which is the probability in $(a)$ in the case $l<n-1$.
\end{prooof}
Theorem~\ref{thm-dual} now easily follows. Indeed, from the last
proposition we see that $\varphi(X)$ is equal in law to the skeleton chain
$\{\mathcal{X}'_k\}_{0\leq k\leq n}$ started from $N$ particles of unit
mass. By the duality relation~\eqref{eq:dualityXY}, reversing $\varphi(X)$
in time yields the process $\varphi(Y)$, which is the state chain of a
fragmentation process.


 \section{Random binary forest representation}
\label{binaryforests}
In this section, we give a second construction of the skeleton chain of the
ternary coalescent in terms of random binary forests. The connection
between random forests and coalescent processes was first observed by
Pitman in~\cite{PIT:99}. In our description, we are guided by Chapter 5.2.3 of
Bertoin~\cite{BER:1}. 
\subsection{Basic definitions on graphs}
We first collect some basic notions on graphs which
will useful for our purpose.

A (undirected) {\it graph} is a pair $G=(V,E)$, where $V$ is a finite set
and $E\subset \{U\subset V : |U|=2\}$. The elements of $V$ are called {\it
  vertices}, the elements of $E$ {\it edges}. The {\it size} of a graph is
the number of vertices $|V|$.  A {\it subgraph} of a graph
$G=(V,E)$ is a graph $H=(V',E')$ with $V'\subset V$ and $E'\subset E$.

Now let $G=(V,E)$ be a graph. Two vertices $v,w$ are {\it adjacent}, if
$\{v,w\}\in E$. The {\it degree} of a vertex $v$ is the number of vertices
adjacent to $v$. 
A sequence $(v_1,e_1,v_2,\ldots,v_m,e_m,v_{m+1})$ such that $m\geq 0$
$v_i\neq v_j$ for $i\neq j$ 
and $e_i = \{v_i,v_{i+1}\}\in E$ for $1\leq
i\leq m$ is called a {\it path}, or also a $v_1$-$v_{m+1}$-path. A {\it
  cycle} is a sequence $(v_1,e_1,\ldots,v_m,e_m,v_1)$ such that $m\geq 2$,
$(v_1,e_1,\ldots,v_{m-1},e_{m-1},v_m)$ is a path and
$e_{m}=\{v_m,v_1\}\in E$. 
We say that two vertices $v,w$ are {\it connected}, if
there exists a $v$-$w$-path. If there is a $v$-$w$-path for any $v,w\in V$,
we say that the graph $G$ is connected. The maximal connected subgraphs of
$G$ are its {\it connected components}. A connected graph without a cycle (as a
subgraph) is called a {\it tree}. In a tree, a {\it leaf} is a vertex of
degree equals $1$, while the vertices of degree greater than $1$ are called
{\it internal} vertices.

We are interested in a special family of trees. A {\it binary tree} is
either a tree consisting of a single vertex only, called the {\it root} of
the tree, or a tree where exactly one vertex has degree $2$, which we then
call the root of the tree, and all the other vertices have degree $3$ or
they are leaves. 
The {\it height} of a vertex $v$ in a binary tree is the number of edges of
the (unique) $v$-$r$-path, where $r$ is the root of the tree. If $v$ is not
a leaf, then there are exactly two vertices $w, w'$ adjacent to $v$ with
height strictly bigger than that of $v$, the {\it children} of $v$. We call
the pair $\{\{v,w\},\{v,w'\}\}$ the {\it outgoing edges} (from
$v$). Finally, a {\it binary forest} is a graph such that its connected
components are binary trees. The leaves or internal vertices of such a
forest are then all those of its tree components.

Observe that a binary forest on $N$ vertices with $m$ tree components has
$N-m$ edges, $(N+m)/2$ leaves, and $(N-m)/2$ internal vertices.

\begin{remark}
 In the literature, a binary tree in our sense is often called a
  (rooted) {\it full labeled} binary tree. The term ``full'' reflects
  the fact that every vertex other than the leaves has two children, and
  ``labeled'' stresses that the vertices are distinguishable. However, we
  will use the term ``labeled'' to indicate a labeling of internal
  vertices.
\end{remark}

\subsection{Dynamics}
Our concern here is to describe the dynamics on the space of binary
forests, which will lead to another representation of the ternary
coalescent.
 
As before let $N=2n+1$. We consider $V=\{1,2,\ldots,N\}$ as a set of
vertices. Given a binary forest on $V$, we enumerate its tree components
according to the increasing order of their roots.
 
We will assign additional labels to all internal vertices of such a
forest. A {\it labeling} of a binary forest on $V$ with $m$ tree
components is a bijective map from the set of $(N-m)/2$ internal vertices
into $\{1,\ldots,(N-m)/2\}$. A labeled binary forest on $V$ is then a
binary forest together with a labeling. Note that internal vertices are
double-labeled, by $V$ and by the labeling
just described.
The set of all labeled binary
forests on $V$ with $m$ tree components is denoted by $\pF(m,N)$. Clearly, 
$\pF(m,N)$ is empty if $m$ is an even number. 
 
For every $1\leq k \leq n$, we define a map $R:\pF(2k-1,N) \longrightarrow
\pF(2k+1,N)$ as follows. For each $\tau\in\pF(2k-1,N)$, select the internal
vertex with the highest label and delete both outgoing edges (and the
label, since the vertex is now a leaf). We obtain a labeled binary forest
with $2k+1$ trees, which we denote by $R(\tau)$.

As the reader might already guess, the map $R$ will be the building block of
the fragmentation mechanism - it breaks the tree with the highest label
into three (new) trees. The reverse dynamic will correspond to the
coagulation mechanism: Out of a binary forest with at least three trees,
pick one leaf and connect it by adding edges to two distinct roots from
other tree components. Then, three trees have merged into one (new) tree,
and the selected leaf has become an internal vertex. Before underlying this
procedure with randomness, let us analyze the map $R$ in detail.
\begin{lemma}
\label{mapR}
For every $1\leq k\leq n$, the map $R:\pF(2k-1,N) \longrightarrow
\pF(2k+1,N)$ is surjective. More precisely, for every
$\tau\in\pF(2k+1,N)$,
\begin{equation*}
\left|\{\tilde{\tau}\in\pF(2k-1,N) : R(\tilde{\tau}) = \tau\}\right| = (n+k+1)k(2k-1).
\end{equation*}
\end{lemma}
\begin{prooof}
  Let $\tau\in\pF(2k+1,N)$. In order to construct a generic
  $\tilde{\tau}\in R^{-1}(\tau)$, pick a leaf $i$ from $\tau$. Write
  $\rho(i)$ for the root of the tree component containing $i$.  Then select
  two roots $j\neq j'$ different from $\rho(i)$, add the edges $\{i,j\}$,
  $\{i,j'\}$ and label the vertex $i$ with the number $n-k+1$. Out of three
  components, we have obtained a new labeled binary tree with root
  $\rho(i)$, which is part of a forest with $2k-1$ trees. Clearly, this
  forest is contained in $R^{-1}(\tau)$. Moreover, different choices of
  $i,j,j'$ give rise to different forests. To finish the proof, note that
  there are $n+k+1$ possible choices for a leaf $i$, and $2k(2k-1)/2$
  possible choices for distinct roots $\{j,j'\}$.
\end{prooof}

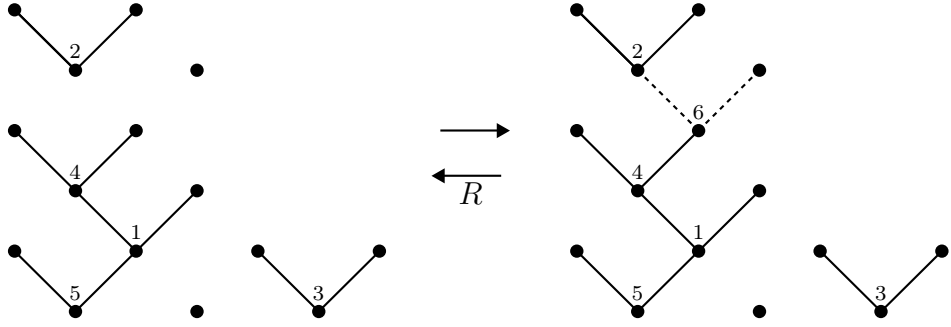
\begin{figure}
\centering
\subfigure{
\psset{xunit=0.8cm,yunit=0.8cm,algebraic=true,dotstyle=o,dotsize=3pt 0,linewidth=0.8pt,arrowsize=3pt 2,arrowinset=0.25}
\begin{pspicture*}(1.6,-0.5)(10.5,5.5)
\psline(3,0)(2,1)
\psline(3,0)(4,1)
\psline(4,1)(3,2)
\psline(3,2)(2,3)
\psline(3,2)(4,3)
\psline(4,1)(5,2)
\psline(7,0)(6,1)
\psline(7,0)(8,1)
\psline(2,5)(3,4)
\psline(4,5)(3,4)
\psline(3,4)(2,5)
\psline(9,3)(10,3)
\psline(9,2.25)(10,2.25)
\rput[t](9.5,2.15){$R$}
\begin{scriptsize}
\psdots[dotsize=5pt 0,dotstyle=*](3,0)
\rput[b](3,0.2){5}
\psdots[dotsize=5pt 0,dotstyle=*](2,1)
\psdots[dotsize=5pt 0,dotstyle=*](4,1)
\rput[b](4,1.2){1}
\psdots[dotsize=5pt 0,dotstyle=*](3,2)
\rput[b](3,2.2){4}
\psdots[dotsize=5pt 0,dotstyle=*](5,2)
\psdots[dotsize=5pt 0,dotstyle=*](2,3)
\psdots[dotsize=5pt 0,dotstyle=*](4,3)
\psdots[dotsize=5pt 0,dotstyle=*](5,0)
\psdots[dotsize=5pt 0,dotstyle=*](7,0)
\rput[b](7,0.2){3}
\psdots[dotsize=5pt 0,dotstyle=*](6,1)
\psdots[dotsize=5pt 0,dotstyle=*](8,1)
\psdots[dotsize=5pt 0,dotstyle=*](3,4)
\rput[b](3,4.2){2}
\psdots[dotsize=5pt 0,dotstyle=*](2,5)
\psdots[dotsize=5pt 0,dotstyle=*](4,5)
\psdots[dotsize=5pt 0,dotstyle=*](5,4)
\psdots[dotsize=5pt 0,dotstyle=triangle*,dotangle=270](10,3)
\psdots[dotsize=5pt 0,dotstyle=triangle*,dotangle=90](9,2.25)
\end{scriptsize}
\end{pspicture*}}
\subfigure{
\psset{xunit=0.8cm,yunit=0.8cm,algebraic=true,dotstyle=o,dotsize=3pt 0,linewidth=0.8pt,arrowsize=3pt 2,arrowinset=0.25}
\begin{pspicture*}(1.6,-0.5)(8.5,5.5)
\psline(3,0)(2,1)
\psline(3,0)(4,1)
\psline(4,1)(3,2)
\psline(3,2)(2,3)
\psline(3,2)(4,3)
\psline(4,1)(5,2)
\psline(7,0)(6,1)
\psline(7,0)(8,1)
\psline(2,5)(3,4)
\psline(4,5)(3,4)
\psline(3,4)(2,5)
\psline[linestyle=dashed,dash=2pt 2pt](3,4)(4,3)
\psline[linestyle=dashed,dash=2pt 2pt](4,3)(5,4)
\begin{scriptsize}
\psdots[dotsize=5pt 0,dotstyle=*](3,0)
\rput[b](3,0.2){5}
\psdots[dotsize=5pt 0,dotstyle=*](2,1)
\psdots[dotsize=5pt 0,dotstyle=*](4,1)
\rput[b](4,1.2){1}
\psdots[dotsize=5pt 0,dotstyle=*](3,2)
\rput[b](3,2.2){4}
\psdots[dotsize=5pt 0,dotstyle=*](5,2)
\psdots[dotsize=5pt 0,dotstyle=*](2,3)
\psdots[dotsize=5pt 0,dotstyle=*](4,3)
\rput[b](4,3.2){6}
\psdots[dotsize=5pt 0,dotstyle=*](5,0)
\psdots[dotsize=5pt 0,dotstyle=*](7,0)
\rput[b](7,0.2){3}
\psdots[dotsize=5pt 0,dotstyle=*](6,1)
\psdots[dotsize=5pt 0,dotstyle=*](8,1)
\psdots[dotsize=5pt 0,dotstyle=*](3,4)
\rput[b](3,4.2){2}
\psdots[dotsize=5pt 0,dotstyle=*](2,5)
\psdots[dotsize=5pt 0,dotstyle=*](4,5)
\psdots[dotsize=5pt 0,dotstyle=*](5,4)
\end{scriptsize}
\end{pspicture*}}
\caption{From left to right (right to left) one step in the coagulation
  (fragmentation) mechanism is shown. For simplicity, only the labeling of
  the internal vertices is depicted. On the left, the leaf with number $6$
  on the right is chosen as well as two roots from different tree
  components. They are connected by two edges visualized by the dashed lines on
the right side.}
\end{figure}

\begin{remark}
\label{forestsizes}
  Applying the map $R$ at most $n$ times destructs a labeled binary forest
  into its single vertices. Due to the recursive structure of trees,
  this method enables one to compute various combinatorial quantities. For
  example, using $|\pF(N,N)| = 1$ and iteratively the identity
\begin{equation*}
\left|\pF(2k-1,N)\right| = (n+k+1)k(2k-1)\left|\pF(2k+1,N)\right|
\end{equation*}
provided by Lemma~\ref{mapR}, one obtains for $k=2,\ldots, n+1$
\begin{equation*}
  \left|\pF(2k-1,N)\right| = \frac{2^{k-(n+1)}n\,(2n+1)!\,(2n-1)!\,(k-2)!}{(n+k)!\,(k-1)!\,(2k-3)!}.
\end{equation*}
In the case $k=1$,
\begin{equation*}
  \left|\pF(1,N)\right| =\frac{2^{-n}(2n)!(2n+1)!}{(n+1)!}= 2^{-n}(2n+1)!\,n!\,C_n,
\end{equation*}
where $C_n = (2n)!/((n+1)!\,n!)$ is the $n$th {\it Catalan} number. Since
there are $n!$ different labelings of internal vertices, we deduce that the
number of binary trees on $V$ is given by $2^{-n}(2n+1)!\ C_n$. 
\end{remark}

\subsection{From forests to  mass partitions}
Denote by $R^k$ the $k$th concatenation of $R$, where $R^0$ is the
identity map. We randomize the input by endowing the space $\pF(1,N)$ with the uniform
probability measure and interpret the maps $R^k$ as random variables
\begin{equation*}
R^k:\pF(1,N)\longrightarrow \pF(2k+1,N), \quad k=0,\ldots,n.
\end{equation*}
In words, $R^k(\tau)$ is the forest with $2k+1$ tree components which
arises from $\tau\in\pF(1,N)$ by picking the $k$ internal vertices with the
highest labels and deleting their outgoing edges. By induction, we deduce
from Lemma~\ref{mapR} that $R^k$ obeys the uniform law on the space
$\pF(2k+1)$, for each $k$.  We then consider the random variables
\begin{equation*}
|R^k|^\downarrow:\pF(1,N)\longrightarrow \ump_{2n+1}, \quad k=0,\ldots,n,
\end{equation*}
where for a tree $\tau\in\pF(1,N)$, $|R^k|^\downarrow(\tau) = {\bf
  s}=(s_1,\ldots,s_{2k+1})\in\ump_{2n+1}$ is the sequence of the sizes of
the tree components, ranked in decreasing order.

Turning back to the ternary coalescent, let $\mathcal{X}'_k$,
$k=0,\ldots,n$, denote the skeleton chain started
from $N$ particles of unit mass. Its connection to the sizes of the tree
components is given by
\begin{proposition}
\label{randtrees-prop} 
The sequence of random variables $\{|R^{n-k}|^\downarrow\}_{0\leq k\leq n}$ is
the state chain of the ternary coalescent, that is
\begin{equation*}
\left(|R^n|^\downarrow,|R^{n-1}|^\downarrow,\ldots,|R^0|^\downarrow\right)
 \eql \left(\mathcal{X}'_0,\ldots,\mathcal{X}'_n\right).
\end{equation*}
\end{proposition}
\begin{prooof}
For each tree $\tau\in\pF(1,N)$, the forest $R^n(\tau)$ has no edges, so $|R^n|^\downarrow=(1,\ldots,1) = \mathcal{X}'_0$. 
Note that given $|R^{l}|^\downarrow = {\bf s}=(s_1,\ldots,s_{2l+1})$ for
some $1\leq l\leq n $, the mass partition $|R^{l-1}|^\downarrow$ is
obtained from ${\bf s}$ by replacing three elements $s_i$,$s_j$,$s_k$,
where $i$,$j$,$k$ are pairwise distinct, by their sum. Furthermore, observe
that the random variables $R^k$, $l\leq k\leq n$, are measurable with
respect to the sigma-field generated by $R^l$. In particular, by
Proposition~\ref{coal-proc-prop}, the claim follows if we show that for
every $0\leq l <n$, for every ${\bf
  s}=(s_1,\ldots,s_{2(n-l)+1})\in\ump_{2n+1}$ and for every $3$-subset
$\{i,j,k\}\subset\{1,\ldots,2(n-l)+1\}$,
\begin{equation*}
\pP\left(|R^{n-l-1}|^\downarrow= s^{i\oplus j\oplus
    k}\ |\ R^{n-l},\, |R^{n-l}|^\downarrow = {\bf s} \right) =
\frac{s_i+s_j+s_k+3}{(2n+1-l)2(n-l)(2(n-l)-1)}.
\end{equation*}
Take a forest $\tau\in\pF(2(n-l)+1,N)$. We work conditionally on
$R^{n-l}=\tau$. By our observation above, $R^{n-l-1}$ is uniformly
distributed on the set of $(2n+1-l)(n-l)(2(n-l)-1)$ forests which can be
obtained from $\tau$ in the way described in Lemma~\ref{mapR}. We
write $\tau_1,\ldots,\tau_{2(n-l)+1}$ for the tree
components of $\tau$. For every $3$-subset $\{a,b,c\}\subset
\{1,\ldots,2(n-l)+1\}$, the probability that the leaf $i$ is picked in
$\tau_a$ and the roots are chosen from $\tau_b$ and
$\tau_c$ is therefore
\begin{equation*}
\frac{|\tau_a|+1}{2}\times \frac{1}{(2n+1-l)(n-l)(2(n-l)-1)}.
\end{equation*}
Hence the probability that $R^{n-l-1}$ evolves from $\tau$ by
merging the trees $\tau_a$, $\tau_b$ and
$\tau_c$, that is the probability that the leaf $i$ is picked in
either $\tau_a$, $\tau_b$ or $\tau_c$ and connected to the roots
of the other two components is
\begin{equation*}
 \frac{|\tau_a|+|\tau_b|+|\tau_c|+3}{(2n+1-l)2(n-l)(2(n-l)-1)}.
\end{equation*}
\end{prooof}
As a consequence of Theorem~\ref{thm-dual}, the time-reversed process
$\{|R^{k}|^\downarrow\}_{0\leq k\leq n}$ is a fragmentation chain with dislocation law $\mu$.

\begin{remark}
\label{1dstat-2}
Adapting the proof of Corollary 5.7 in~\cite{BER:1} to our situation, we
find another  way to prove Corollary~\ref{1dstat-cor}, based
on the binary forest representation. 
Namely, with $m=2(n-l)+1$ and ${\bf s} = (s_1,\ldots, s_m)\in\ump_{2n+1}$, there are 
\begin{equation*}
\frac{1}{m!}{2n+1\choose s_1,\ldots,s_m} = \frac{(2n+1)!}{m!\,s_1!\cdots s_m!}
\end{equation*}
possibilities to partition the set of vertices $\{1,\ldots,2n+1\}$ into non-empty
disjoint sets $E_i$, $i=1,\ldots,m$, such that $|E_i| = s_i$ and
$\min E_i< \min E_j$ for $i<j\leq m$.  Without
labeling internal vertices, the number of binary tree structures
which can be attached to $E_i$ is $|\pF(1,s_i)|/((s_i-1)/2)!$.  Having
chosen a binary tree structure for each $E_i$, there are $l!$ possible ways
to label the $l$ internal vertices. Recall that the tree components
of a forest are enumerated in increasing order of their roots. It follows
that the number of binary forests $\tau\in\pF(m,2n+1)$ with tree
components $\tau_i$ such that $|\tau_i|=s_i$ is given by
\begin{equation*}
  \frac{(2n+1)!\,l!}{m!}\,\prod_{i=1}^m\frac{|\pF(1,s_i)|}{s_i!\left(\frac{s_i-1}{2}\right)!}.
\end{equation*}
Since $R^{n-l}$ is uniformly distributed on $\pF(m,2n+1)$, we
deduce from Proposition~\ref{randtrees-prop} that
\begin{equation*}
  \pP\left(\mathcal{X}'_l = (s_1,\ldots,s_m)\right) = \frac{\gamma({\bf s})}{|\pF(m,2n+1)|}
  \frac{(2n+1)!\,l!}{m!}\,\prod_{i=1}^m\frac{|\pF(1,s_i)|}{s_i!\left(\frac{s_i-1}{2}\right)!},
\end{equation*}
where $\gamma({\bf s})$ has been defined in Section~\ref{tc-mono}.
Plugging in the values for $|\pF(m,2n+1)|$ and $|\pF(1,s_i)|$ from
Remark~\ref{forestsizes} results in the expression obtained in
Corollary~\ref{1dstat-cor}.
\end{remark}

\subsection{Encoding forests by paths}
\begin{figure}
\centering
\subfigure{
\psset{xunit=0.8cm,yunit=0.8cm,algebraic=true,dotstyle=o,dotsize=3pt 0,linewidth=0.8pt,arrowsize=3pt 2,arrowinset=0.25}
\begin{pspicture*}(1.5,-0.5)(9,4)
\psline(3,0)(2,1)
\psline(3,0)(4,1)
\psline(4,1)(3,2)
\psline(3,2)(2,3)
\psline(3,2)(4,3)
\psline(4,1)(5,2)
\psline(7,0)(6,1)
\psline(7,0)(8,1)
\begin{scriptsize}
  \psdots[dotsize=5pt 0,dotstyle=*](3,0) \rput[b](3,0.2){0}
  \psdots[dotsize=5pt 0,dotstyle=*](2,1) \rput[b](2,1.2){1}
  \psdots[dotsize=5pt 0,dotstyle=*](4,1) \rput[b](4,1.2){2}
  \psdots[dotsize=5pt 0,dotstyle=*](3,2) \rput[b](3,2.2){3}
  \psdots[dotsize=5pt 0,dotstyle=*](5,2) \rput[b](5,2.2){6}
  \psdots[dotsize=5pt 0,dotstyle=*](2,3) \rput[b](2,3.2){4}
  \psdots[dotsize=5pt 0,dotstyle=*](4,3) \rput[b](4,3.2){5}
  \psdots[dotsize=5pt 0,dotstyle=*](5,0) \rput[b](5,0.2){7}
  \psdots[dotsize=5pt 0,dotstyle=*](7,0) \rput[b](7,0.2){8}
  \psdots[dotsize=5pt 0,dotstyle=*](6,1) \rput[b](6,1.2){9}
  \psdots[dotsize=5pt 0,dotstyle=*](8,1) \rput[b](8,1.2){10}
\end{scriptsize}
\end{pspicture*}}
\subfigure{
\psset{xunit=0.7cm,yunit=0.7cm,algebraic=true,dotstyle=o,dotsize=3pt 0,linewidth=0.8pt,arrowsize=3pt 2,arrowinset=0.25}
\begin{pspicture*}(-0.5,-3.5)(11.5,2.5)
\psaxes[labelFontSize=\scriptscriptstyle,xAxis=true,yAxis=false,Dx=1,Dy=1,ticksize=-2pt 0,subticks=2]{-}(0,0)(-0.5,-3.5)(11.5,2.5)
\psline(0,0)(1,1)
\psline(1,1)(2,0)
\psline(2,0)(4,2)
\psline(4,2)(8,-2)
\psline(8,-2)(9,-1)
\psline(9,-1)(11,-3)
\begin{scriptsize}
  \psdots[dotsize=7pt 0,dotstyle=x](0,0) \psdots[dotsize=7pt
  0,dotstyle=x](7,0) \psdots[dotsize=7pt 0,dotstyle=x](8,0)
  \rput[t](11,-3.1){(11,-3)}
\end{scriptsize}
\end{pspicture*}}
\caption{On the left side a binary plane forest on $11$ vertices with $3$
  tree components is shown, where the vertices are enumerated by a depth-first
  search. The corresponding Lukasiewicz path is depicted on the right
  side. The crosses indicate where the coding of a new tree
  starts.}
\label{fig:PBF-LP}
\end{figure}
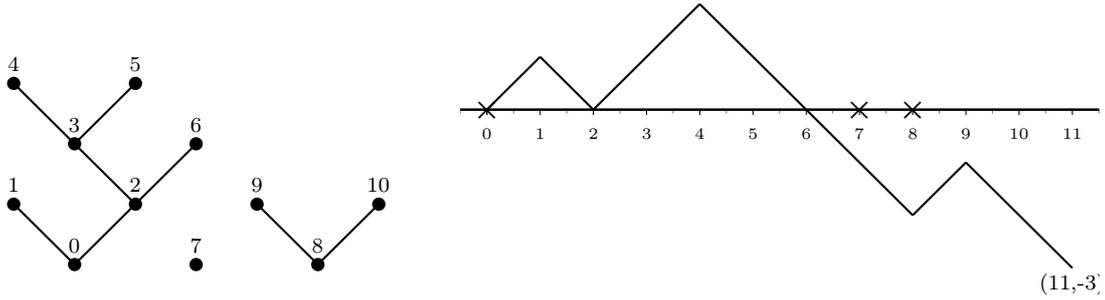

We conclude our discussion of binary forests by illustrating a direct
connection to the random walk representation. 
Here, it is more convenient to consider (rooted unlabeled)
plane trees and forests. In a plane forest vertices are regarded as
indistinguishable, but the set of children for each vertex is ordered, as
well as the set of roots of the different tree components. The ordering induces
serveral natural enumerations of the vertices. For example, one of them is
provided by the order in which the vertices are visited by a depth-first
search, see Figure~\ref{fig:PBF-LP}. More on this can be found in Chapter 6.2 of Pitman~\cite{PIT:StF}.

We will look at (full) binary plane forests. To relate them to the
binary forests considered above, note that the number of binary
plane forests on $N$ vertices with $k$ tree components is equal to
\begin{equation*}
 \frac{2^{(N-k)/2}k!}{N!((N-k)/2)!}|\pF(k,N)|,
\end{equation*}
since there are $2^{(N-k)/2}$ possible orderings of the children of the
internal vertices of a forest in $\pF(k,N)$, $k!$ orderings of the roots,
but neither vertices are labeled nor there is an additional identification
of internal vertices.  Clearly the ternary coalescent with a monodisperse
initial configuration can also be realized on the space of binary plane
forests, with the same dynamics.

There are various possibilities to code plane trees and forests by discrete
functions. For a (finite) plane tree $\theta$ on $N$ vertices, one common way is to look  
at its Lukasiewicz path $\{x_l\}_{0\leq l\leq N}$. Denoting by
$v_0,\ldots,v_{N-1}$ the vertices of $\theta$ listed in the order of a
depth-first search and by $k(v)$ the number of children of vertex $v$, one
defines
\begin{equation*}
x_j = \sum_{i=0}^{j-1}(k(v_i)-1),\quad 0\leq j\leq N.
\end{equation*}
Note that $x_0=0$, $x_N=-1$, and
\begin{equation}
\label{eq:lukpath}
x_j-x_{j-1} = k(v_{j-1})-1,\quad 1\leq j\leq N.
\end{equation}
It is easy to see that there is a bijection between Lukasiewicz paths and
rooted plane trees. A sequence of such trees may then by encoded by gluing
together the corresponding Lukasiewicz paths, retaining the
relationship~\eqref{eq:lukpath}. In other words, the coding of the next
tree starts if a new minimum is attained.

Turning to random trees, it follows from Proposition 1.4 of Le
Gall~\cite{LG} that a Galton-Watson tree with offspring distribution
$\eta(k) = 1/2 (\delta_0(k)+\delta_2(k))$, conditioned to have total
progeny size $N$, is distributed according to a tree chosen uniformly at
random among the set of all binary plane trees on $N$ vertices.
Further, the corresponding Lukasiewicz path
tree is distributed as the path of simple random walk on $\mathbb{Z}$ up to
time $N$, conditioned on $H_{-1}=N$ (see Corollary 1.6 of~\cite{LG}).

We then realize that for an integer $0\leq
l\leq n$, the path of simple random walk up to time $N$, conditioned on
$H_{-(2l+1)}=N$, encodes a forest distributed uniformly over all binary
plane forests on $N$ vertices with $2l+1$ tree components. In particular,
the sequence of the sizes of the tree components is
distributed as the sequence of the ladder epochs of the conditioned random
walk path, if both are put in random uniform order, say. However, the
sequence of coding functions induced by the above dynamics on the space of binary
forests is not directly related to the sequence of paths
of the random walk representation. In this sense, the connection between
the two representations is only static.
\section{Asymptotics of the ternary coalescent}
Having concrete realizations at hand, we are now able to investigate
asymptotic properties of the ternary coalescent process. Let us write
$\mathcal{X}^{[N]} = (\mathcal{X}^{[N]}(t),\,t\geq 0)$ for the coalescent
with kernel $\kappa$ started from the monodisperse configuration
$(1,\ldots,1)$ consisting of $N=2n+1$ atoms of unit mass, and put
$\mathcal{X}_k'^{[N]} = \mathcal{X}^{[N]}(T_k)$, $k=0,\ldots,n$.  The
number of particles at time $t\geq 0$ is denoted by $\#^{[N]}(t)$, and the
number of jumps up to time $t$ by $J^{[N]}(t)$.

We will consider the space of mass partitions with
total mass bounded by $1$,
\begin{equation*}
\massp = \left\{{\bf s}=(s_1,s_2,\ldots) : s_1\geq s_2\geq\ldots\geq
0, \sum_{i=1}^\infty s_i \leq 1\right\},
\end{equation*}
and the subset $\mathbb{S}_1 \subset \massp$ of sequences with
$\sum_{i=1}^\infty s_i=1$.  We equip $\massp$ with the uniform
distance. The induced topology coincides with that of pointwise convergence
and turns $\massp$ into a compact space. The $l_1$-distance induces a finer
topology. However, if $({\bf s}_n,\, n\in\mathbb{N})$ is some sequence in
$\massp$ converging pointwise to ${\bf s}\in\mathbb{S}_1$, then
the convergence does also hold in the $l_1$-sense, as it can be easily
deduced from Scheff\'e's lemma.  Therefore, on $\mathbb{S}_1$ all these
types of convergence are equivalent.

We turn to our main result of this section. Recall that the standard
additive coalescent $\mathfrak{X}=(\mathfrak{X}(t),\, t\in\mathbb{R})$ is the
unique additive coalescent process such that for each $t\in\mathbb{R}$,
$\mathfrak{X}(t)$ has the law of the ranked sequence ${\tt a}_1\geq {\tt
  a}_2\geq ...$ of the atoms of a Poisson random measure on $(0,\infty)$
with intensity measure $\Lambda(da)=e^{-t}da/\sqrt{2\pi a^3}$, conditioned
on $\sum_{i=1}^\infty {\tt a_i} = 1$. We refer to~\cite{ALDPIT}
and~\cite{EVPIT} for background.

\begin{theorem}
\label{convtoaddcoal-thm}
As $n\rightarrow\infty$, the $\mathbb{S}_1$-valued process
\begin{equation*}
t\mapsto \frac{1}{N}\mathcal{X}^{[N]}(e^t/N^{3/2}),\quad t\in\mathbb{R},
\end{equation*}
converges in the sense of finite-dimensional distributions towards the
standard additive coalescent.
\end{theorem}

Here, the multiplication with $1/N$ is meant element-wise.
At first glance the convergence may look surprising, since the standard
additive coalescent is a {\it binary} coalescent that arises
as a limit of additive coalescent processes as follows (Evans and Pitman~\cite{EVPIT}). Let
$\mathfrak{X}^{[n]} = (\mathfrak{X}^{[n]}(t),\,t\geq 0)$ be the stochastic
coalescent with binary coagulation kernel
\begin{equation*}
\tilde{\kappa}(r,s) = r+s,\quad r,s > 0,
\end{equation*}
started from the monodisperse configuration with $n$ atoms, each of mass
$1/n$. Then, as $n\rightarrow\infty$, the time-shifted processes
$(\mathfrak{X}^{[n]}(t+(1/2)\ln n),t\geq -(1/2)\ln n)$ converge weakly to
$\mathfrak{X}$.  

However, our convergence result concerns only the finite-dimensional
laws. For the one-dimensional distributions, one might expect a result in
this direction if one compares the one-dimensional statistics of the
skeleton chains of the ternary and the additive coalescent
$\mathfrak{X}^{[n]}$. The states of the additive coalescent can be
expressed in terms of independent standard Borel variables (see for example
$(30)$ in~\cite{EVPIT}), which have a similar tail behavior as the hitting
time $H_k$. For the finite-dimensional laws, an analysis of the first
hitting time distribution shows that a ``true'' ternary coagulation step,
i.e. the event that three particles merge which are all of a size
comparable to $n$, only occurs with negligible probability. Therefore,
under the rescaling, the process looks more like a binary coalescent.

Let us briefly comment on the scaling in the theorem. To obtain a limit for
the normalized sequence of masses, the number of atoms must be of order
$\sqrt{n}$. We refer to Lemma~\ref{hittingtimes-lemma2} for
a better understanding. As Lemma~\ref{ternarycoal-asymplemma} shows, if the
process $\mathcal{X}^{[N]}$ runs for time $t/N^{3/2}$, then the amount of
particles has typically reduced from $N$ to about $\sqrt{N}/t$. Note that
when approximating the standard additive coalescent with the processes
$\mathfrak{X}^{[n]}$ starting from $n$ atoms of mass $1/n$, the macroscopic
picture appears at times $t+(1/2)\ln n$, at which there are about
$\sqrt{n}/e^{t}$ particles. Here, roughly speaking, the standard Borel law
plays the role of the hitting time distribution. Precise statements can be
found in the books of Pitman~\cite{PIT:StF}, Chapter 10.3, and 
Bertoin~\cite{BER:1}, Chapter 5.3.

We shall present three different ways to obtain convergence for the rescaled
ternary coalescent of which we discuss two in detail. The first more
general method will lead to one-dimensional convergence in
Proposition~\ref{prop-onedimconv}. It relies on the observation that the
distribution of the hitting time $H_k$ is in the domain of attraction of a
stable$(1/2)$ law. Then a size-biased reordering is used to construct the
limiting mass partition.  The second method resulting in finite-dimensional
convergence (and therefore in the proof of the theorem) is more specialized
to our situation. It is based on the identification of configurations with
mass partitions via paths, as described in
Section~\ref{mc-and-fcts}.
Since the two methods do not rely on each other, the reader in a hurry may
safely skip Section~\ref{poissonmeasures}. In a closing remark we
outline a possible third way to establish finite-dimensional convergence,
using the random binary forest representation.

\subsection{Number of particles}
In order to relate the behavior of $\mathcal{X}^{[N]}$ to that of its
skeleton chain, we prove a limit theorem for the number of particles. As just
remarked, it will become clear later why we choose the spatial scale factor
$N^{-1/2}$.
\begin{lemma}
\label{ternarycoal-asymplemma}
For every $t> 0$, as $n\rightarrow \infty$,
\begin{equation*}
\frac{\#^{[N]}(t/N^{3/2})}{\sqrt{N}}\rightarrow
\frac{1}{t}\quad\mbox{in probability}.
\end{equation*}
\end{lemma}
\begin{prooof}
Using the relation $\#^{[N]}(\cdot)= N - 2J^{[N]}(\cdot)$, the
claim will follow once we show that 
\begin{equation}
\label{eq:asymplemma-1}
\frac{J^{[N]}(t/N^{3/2})}{\sqrt{N}}-\frac{(\sqrt{N}-t^{-1})}{2}\rightarrow
  0\quad\mbox{in probability}.
\end{equation}
Remember that $J^{[N]}(t/N^{3/2}) = \max\{k\in\mathbb{Z}_+ : N^{3/2}T_k\leq
t\}$, where $T_k$ is the $k$th coagulation time given by $T_k \eql
\sum_{i=1}^{k}\alpha(i)^{-1}{\bf e}_i$, the rates $\alpha(i) = \alpha(i,N)$
are as in~\eqref{eq:tc-monorates} and ${\bf e}_1,{\bf e}_2,\ldots$ is a
sequence of independent standard exponential variables. Heuristically,
replacing $T_k$ by its expectation $\sum_{i=1}^k\alpha(i)^{-1}$, the number
of jumps $J^{[N]}(t/N^{3/2})$ should roughly behave as the maximal $k$ such that
$N^{3/2}\sum_{i=1}^k\alpha(i)^{-1} \leq t$. We will show that with the
choice $k_n=n - t^{-1}\sqrt{N}/2$,
\begin{equation}
\label{eq:asymplemma-2}
N^{3/2}\sum_{i=1}^{k_n}\alpha(i)^{-1} = t + o(1),
\end{equation}
where we agree that the sum runs from $1$ to the largest integer below
$k_n$.  First note that
\begin{eqnarray*}
  N^{3/2}\sum_{i=1}^{k_n}\alpha(i)^{-1} &=& N^{3/2}\sum_{i=1}^{k_n}\frac{1}{(N+1-i)(N+1-2i)(N-2i)}\\
  &=& N^{3/2} \left(\sum_{i=1}^{k_n}\frac{1}{(N-i)(N-2i)^2}\right)+O(n^{-1/2}).
\end{eqnarray*}
Furthermore, some simple computations show that for each $\varepsilon > 0$,
\begin{eqnarray*}
  \sum_{i=1}^{k_n}\frac{1}{(N-i)(N-2i)^2}& =&
  \int_0^{k_n}\frac{\dt x}{(N-x)(N-2x)^2} +O(n^{-2})\\ &=&
  \frac{1}{N(N-2k_n)} + O(n^{-2+\varepsilon})\\ &=& \frac{t}{N^{3/2}}
  +O(n^{-2+\varepsilon}).
\end{eqnarray*}
Altogether, we obtain~\eqref{eq:asymplemma-2}. Moreover, since
\begin{equation*}
  \textup{Var}\left(N^{3/2}T_{k_n}\right) =
  N^{3}\sum_{i=1}^{k_n}\alpha(i)^{-2} = O(n^{-1/2}),
\end{equation*}
we deduce that $N^{3/2}T_{k_n}\rightarrow t$ in probability. From
this~\eqref{eq:asymplemma-1} readily follows.
\end{prooof}

\subsection{Mass partitions induced by Poisson measures}
\label{poissonmeasures}
We shall now prove one-dimensional convergence of the ternary coalescent
process. First let us recall some basic facts about mass partitions
obtained from Poisson measures, as provided in Section
2.2.3 of Bertoin~\cite{BER:1}.
Consider a measure $\Lambda$ on $(0,\infty)$ such that
\begin{equation}
\label{eq:pintensity}
\int_0^\infty(1\wedge x)\Lambda(dx) < \infty\quad\mbox{and}\quad
\Lambda((0,\infty)) = \infty.
\end{equation}
Let $M$ be a Poisson random measure on $(0,\infty)$ with intensity
$\Lambda$. From~\eqref{eq:pintensity} it follows that $M$ has almost surely a
countably infinite number of atoms, which we may rank in decreasing order,
\begin{equation*}
{\tt a}_1 \geq {\tt a}_2\geq \ldots> 0.
\end{equation*}
Under condition~\eqref{eq:pintensity}, we further have
\begin{equation*}
\boldsymbol{\varsigma} = \sum_{i=1}^\infty {\tt a}_i < \infty\quad\mbox{almost surely.}
\end{equation*}
In our situation, $\Lambda$ will be non-atomic, which implies that the
atoms ${\tt a}_i$ are almost surely distinct. Furthermore, $\Lambda$ will be of a form that
guarantees the existence of a continuous density of $\boldsymbol{\varsigma}$,
\begin{equation*}
\pP(\boldsymbol{\varsigma} \in dx)=\rho(x)dx,\quad x>0,
\end{equation*}
with $\rho >0$ on $(0,\infty)$.

Given some fixed $x > 0$, we want to transform the atoms $({\tt a}_1,{\tt
a}_2,\ldots)$ into a random mass-partition with total mass $1$ by looking
at $({\tt a}_1/x,{\tt a}_2/x,\ldots)$ conditioned on $\sum_{i=1}^\infty
{\tt a}_i =x$. In order to define the singular conditioning in a proper
way, it is useful to look first at a size-biased reordering $\left({\tt a}^\ast_i,
 i\in\mathbb{N}\right)$ of $({\tt a}_1,{\tt a}_2,\ldots)$. This means that
conditionally on $({\tt a}_1,{\tt a}_2,\ldots)$, we choose an index
$1^{\ast}$ according to
\begin{equation}
\label{eq:size-biased-reordering}
\pP\left(1^\ast = k\,|\,({\tt a}_1,{\tt a}_2,\ldots)\right) = {\tt
  a}_k/{\sum_{i=1}^\infty {\tt a}_i},\quad k\in\mathbb{N},
\end{equation}
set ${\tt a}_1^\ast={\tt a}_{1^\ast}$, remove ${\tt a}_1^{\ast}$ from the sequence and
repeat~\eqref{eq:size-biased-reordering} with this new sequence to obtain
$2^\ast$, set ${\tt a}_2^{\ast}={\tt a}_{2^\ast}$, and so on.
In the following, Propositions 2.3
and 2.4 from~\cite{BER:1} play a major role, so we summarize them for
convenience in the next statement.
\begin{proposition}
\label{prop2.3and2.4BER}
\begin{enumerate}
\item Consider for each $n\in\mathbb{N}\cup\{\infty\}$ a random mass
  partition $S^{(n)}$ with total mass equals one almost surely, and a
  size-biased reordering $S^{(n)\ast}$ of $S^{(n)}$. Then, as $n\rightarrow\infty$, convergence in
  distribution of $S^{(n)}$ to $S^{(\infty)}$  in
  $\massp$ is equivalent to convergence of $S^{(n)\ast}$ to $S^{(\infty)\ast}$
  in the sense of finite-dimensional distributions.
\item In the setting from above, for fixed $x>0$, the conditional law of
  $\left({\tt a}^\ast_1, {\tt a}^\ast_2,\ldots\right)$ given
$\boldsymbol{\varsigma}\in[x,x+\varepsilon]$ has a weak limit in the sense
of convergence of finite-dimensional distributions as
$\varepsilon\downarrow 0$, denoted by $\pP_x^\ast$, which is determined
by the following Markov-type property:
\begin{equation*}
\pP_x^{\ast}\left({\tt a}_1^\ast \in dy\right) =
\frac{y\rho(x-y)}{x\rho(x)}\Lambda(dy),\quad 0<y<x,
\end{equation*}
and the conditional distribution of $({\tt a}_2^\ast,{\tt a}_3^\ast,\ldots)$ under
$\pP_x^\ast$ given ${\tt a}_1^\ast=y$ is $\pP_{x-y}^\ast$. Under $\pP_x^\ast$,
$\sum_{i=1}^\infty {\tt a}_i^\ast =x$ almost surely.
\end{enumerate}
\end{proposition}
Having Proposition~\ref{prop2.3and2.4BER} in mind, we call a random sequence
$({\tt a}_1, {\tt a}_2,\ldots)$ that results from the decreasing
rearrangement of $\left({\tt a}^\ast_i, i\in\mathbb{N}\right)$ under
$\pP_x^\ast$ the ranked sequence of the atoms of a Poisson random measure
on $(0,\infty)$ with intensity $\Lambda$, conditioned on
$\sum_{i=1}^\infty{\tt a}_i = x$. We leave it to the reader to check that
$({\tt a}^\ast_1,{\tt a}^\ast_2,\ldots)$ is then a size-biased reordering
of $({\tt a}_1, {\tt a}_2,\ldots)$, in the sense from above. 

Let $\xi_1,\xi_2,\ldots$ be a sequence of
independent copies of $H_{-1}$. Recall that by Lemma~\ref{hittingtimes-lemma}, as
$l\rightarrow\infty$,
\begin{equation}
\label{eq:xi-asymp}
\pP\left(\xi_1 = 2l+1\right) \sim \frac{1}{2}\sqrt{\frac{1}{\pi l^3}}.
\end{equation} 
For $k\in\mathbb{Z}_+$, let
$\Sigma_{2k+1}=\xi_1+\ldots+\xi_{2k+1}$, and denote by $S^{(2k+1,N)}$ a random
mass partition distributed as the rearrangement in decreasing order of
$\xi_1/N,\ldots,\xi_{2k+1}/N$, conditionally on
$\Sigma_{2k+1}=N$. As a special case of Corollary 2.2 in~\cite{BER:1} we
have
\begin{lemma}
\label{hittingtimes-lemma2}
Fix $b > 0$. Then $S^{(2k+1,N)}$ converges in distribution on $\massp$ as $k$,
$n\rightarrow\infty$ with $k\sim bn^{1/2}$ to the ranked sequence $({\tt a}_1,{\tt
  a}_2,\ldots)$ of the atoms of a Poisson random
measure on $(0,\infty)$ with intensity $\Lambda(da)=
b\pi^{-1/2}a^{-3/2}da$, conditioned on $\sum_{i=1}^\infty{\tt a}_i = 1$.
\end{lemma}
\begin{prooof}
For $k\leq n$, denote by $(\xi^\ast_{1,N},\ldots,\xi^\ast_{2k+1,N})$ a
  $(2k+1)$-tuple distributed as a size-biased reordering of
  $(\xi_1,\ldots,\xi_{2k+1})$ given $\Sigma_{2k+1}=N$. 
It easily follows that for $l=0,\ldots,n-k$,
\begin{eqnarray*}
  \pP\left(\xi^{\ast}_{1,N}=2l+1\right) &=&
  \frac{(2k+1)(2l+1)}{N}\pP\left(\xi_1=2l+1\,|\,\Sigma_{2k+1}=N\right)\\
  & =&
  \frac{(2k+1)(2l+1)}{N}\pP\left(\xi_1=2l+1\right)\times\frac{\pP\left(\Sigma_{2k}=N-(2l+1)\right)}{\pP\left(\Sigma_{2k+1}=N\right)}.  
\end{eqnarray*}
If we fix $a\in(0,1)$, $b>0$ and let
$l,k,n$ tend to infinity with $l\sim an$, $k\sim b\,n^{1/2}$, we obtain from~\eqref{eq:xi-asymp} 
\begin{equation}
\label{eq:xi-asymp2}
\frac{(2k+1)(2l+1)}{N}\pP\left(\xi_1 = 2l+1\right) \sim b\,\pi^{-1/2}a^{-1/2}n^{-1}.
\end{equation}
Setting $g_k= 8\pi^{-1}k^2$, we see again by~\eqref{eq:xi-asymp} that for $k\rightarrow\infty$,
\begin{equation*}
(2k+1)\pP\left(\xi_1 > g_k\right) \sim 1.
\end{equation*}
Moreover, since $k\sim b\,n^{1/2}$, we have as $k$, $n\rightarrow\infty$
\begin{equation*}
\frac{g_k}{N}\sim \frac{4b^2}{\pi}.
\end{equation*}
It then follows from the theory of stable laws (see Breiman~\cite{BRE}, Chapters 9 and 14) that 
\begin{equation*}
  \frac{\Sigma_{2k+1}}{N}\rightarrow \boldsymbol{\varsigma}\quad\mbox{in
    distribution as } k, n\rightarrow\infty,
\end{equation*}
where $\boldsymbol{\varsigma}=\sum_{i=1}^\infty{\tt a}_i$ and $({\tt a}_1,{\tt
  a}_2,\ldots)$ is the ranked sequence of the atoms of a Poisson random
measure on $(0,\infty)$ with intensity $b\pi^{-1/2}a^{-3/2}da$. In particular,
$\boldsymbol{\varsigma}$ is stable$(1/2)$ and has a smooth density $\rho(x) =
\pP(\boldsymbol{\varsigma}\in dx)/dx$, which is strictly positive on
$(0,\infty)$. Using $l\sim an$, we infer from
Gnedenko's local limit theorem (see Gnedenko, Kolmogorov~\cite{GK}, p. 236) that
\begin{equation*}
  \frac{\pP\left(\Sigma_{2k} =
      N-(2l+1)\right)}{\pP\left(\Sigma_{2k+1}=N\right)} \sim \frac{\rho(1-a)}{\rho(1)}.
\end{equation*}
Together with~\eqref{eq:xi-asymp2} this shows
\begin{equation*}
  \pP\left(\xi^{\ast}_{1,N} = 2l+1\right)\sim a\frac{\rho(1-a)}{n\rho(1)}b\pi^{-1/2}a^{-3/2}.
\end{equation*}
In particular, $\xi^{\ast}_{1,N}/N$ converges weakly as $n$,
$k\rightarrow\infty$ with $k\sim bn^{1/2}$ towards the law
\begin{equation*}
a\frac{\rho(1-a)}{\rho(1)}b\pi^{-1/2}a^{-3/2}da,\quad a\in (0,1).
\end{equation*}
Now observe that given $\xi^{\ast}_{1,N} = 2l+1$ for some $l=0,\ldots,n-k$,
we have equality in law
\begin{equation*}
\left(\xi^{\ast}_{2,N},\ldots,\xi_{2k+1,N}^{\ast}\right)\eql\left(\xi^{\ast}_{1,N-(2l+1)},\ldots,\xi_{2k,N-(2l+1)}^{\ast}\right).
\end{equation*}
By iterating the argument from above, we may therefore deduce from
the second part of Proposition~\ref{prop2.3and2.4BER} that the limit law of
$\left(\xi^{\ast}_{1,N}/N,\ldots,\xi^{\ast}_{2k+1,N}/N\right)$ in the sense
of finite-dimensional distributions as $n$, $k\rightarrow\infty$ with
$k\sim bn^{1/2}$ is given by the law of a size-biased reordering $({\tt
  a}^\ast_1,{\tt a}^\ast_2,\ldots)$ of $({\tt a}_1,{\tt a}_2,\ldots)$,
where the latter sequence is as in the statement. Also, we have that
$\sum_{i=1}^\infty{\tt a}^\ast_i = 1$ almost surely. From the first part of
the same Proposition it then follows that the (ranked) random mass partition
$S^{(2k+1,N)}$ converges in distribution to the ranked sequence of atoms
$({\tt a}_1,{\tt a}_2,\ldots)$, conditioned on $\boldsymbol{\varsigma} =1$.
\end{prooof}
For the skeleton chain $\mathcal{X}'^{[N]}$, we derive the following consequence.
\begin{corollary}
\label{skeletonconvergence}
  Fix $b>0$. If $n$, $k \rightarrow\infty$ with $k\sim bn^{1/2}$, then
  $(1/N)\mathcal{X}_{n-k}'^{[N]}$ converges in distribution on
  $\mathbb{S}_1$ to the ranked sequence $({\tt a}_1,{\tt a}_2,\ldots)$ of
  the atoms of a Poisson random measure on $(0,\infty)$ with intensity
  $\Lambda(da)= b\pi^{-1/2}a^{-3/2}da$, conditioned on
  $\sum_{i=1}^\infty{\tt a}_i = 1$.
\end{corollary}
\begin{prooof}
This follows from Proposition~\ref{1dstat-prop-RW} together with the
last lemma.
\end{prooof}
Combining the corollary with the weak convergence result for the number of
particles, we easily obtain one-dimensional convergence. 
\begin{proposition}
\label{prop-onedimconv}
Fix $t>0$. Then
\begin{equation*}
\frac{1}{N}\mathcal{X}^{[N]}(t/N^{3/2})
\end{equation*}
converges in distribution on $\mathbb{S}_1$ to the ranked sequence $({\tt
  a}_1, {\tt a}_2,\ldots)$ of the atoms of a Poisson random measure on
$(0,\infty)$ with intensity
\begin{equation*}
\frac{t^{-1}}{\sqrt{2\pi a^3}}da,\quad a>0,
\end{equation*}
conditioned on $\sum_{i=1}^\infty{\tt a}_i = 1$.
In particular, the one-dimensional distributions of the process 
\begin{equation*}
t\mapsto \frac{1}{N}\mathcal{X}^{[N]}(e^t/N^{3/2}),\quad t\in\mathbb{R},
\end{equation*}
converge to those of the standard additive coalescent. 
\end{proposition}
\begin{prooof}
Let $k_n = n-J^{[N]}(t/N^{3/2})$. Then  $\mathcal{X}^{[N]}(t/N^{3/2}) = \mathcal{X}_{n-k_n}'^{[N]}$,
so we may show convergence for $(1/N)\mathcal{X}_{n-k_n}'^{[N]}$. From
Lemma~\ref{ternarycoal-asymplemma} it follows that as $n\rightarrow \infty$,
\begin{equation}
\label{eq:onedimconv-1}
\frac{k_n}{\sqrt{n}} \rightarrow \frac{t^{-1}}{\sqrt{2}}\quad\mbox{in probability}.
\end{equation}
Furthermore, we know from Corollary~\ref{skeletonconvergence} that if $l_n$
is a deterministic sequence of integers with
$l_n\sim\sqrt{n}t^{-1}/\sqrt{2}$, then we have the asserted
convergence for $(1/N)\mathcal{X}_{n-l_n}'^{[N]}$.

It therefore remains to argue that we may replace $l_n$ by the random
sequence $k_n$. To this end, recall that $\massp$ is a compact metric
space, so by Prohorov's theorem (see Billingsley~\cite{BILL}, Section 6)
the space of probability measures on $\massp$ is relatively compact, and we
only have to show convergence on $\massp$ in the sense of
finite-dimensional distributions. Since all our random mass partitions 
lie in $\mathbb{S}_1$ almost surely, this leads to convergence in
distribution on $\mathbb{S}_1$.  Denote by $x_i^{[N]}$ the $i$th component
of $(1/N)\mathcal{X}_{n-k_n}'^{[N]}$. Finite-dimensional convergence on
$\massp$ is equivalent to say that for each $j\in\mathbb{N}$, 
\begin{equation*}
\left(x_1^{[N]},x_1^{[N]}+x_2^{[N]},
\ldots,x_1^{[N]}+\ldots+x_j^{[N]}\right)
\end{equation*}
converges in distribution towards $({\tt a}_1,{\tt a}_1+{\tt
  a}_2,\ldots,{\tt a}_1+\ldots+{\tt a}_j)$, where $({\tt a}_1,{\tt
  a}_2,\ldots)$ is distributed as the rearrangement in decreasing order of $\left({\tt a}^\ast_i, i\in\mathbb{N}\right)$ under
$\pP_1^\ast$, see Proposition~\ref{prop2.3and2.4BER} (i).
This follows if we show that for all $j\in\mathbb{N}$ and $\lambda_i \geq 0$, as
$n\rightarrow \infty$,
\label{eq:onedimconv-2}
\begin{equation}
\pE\left[\exp\left(-\sum_{i=1}^j\lambda_i\left(x_1^{[N]}+\ldots+ x_i^{[N]}\right)\right)\right]\rightarrow
\pE\left[\exp\left(-\sum_{i=1}^j\lambda_i\left({\tt a}_1+\ldots+{\tt a}_i\right)\right)\right].
\end{equation}
Denote by $f:\massp\rightarrow(0,1]$
the function 
\begin{equation*}
f({\bf s})= \exp\left(-\sum_{i=1}^j\lambda_i\left(s_1+\ldots+s_i\right)\right), \quad {\bf s}=(s_1,s_2,\ldots)\in\massp.
\end{equation*}
Note that $f((1/N)\mathcal{X}^{[N]}(t)) \geq f((1/N)\mathcal{X}^{[N]}(s))$
almost surely whenever $t\leq s$. By~\eqref{eq:onedimconv-1} we can find
deterministic sequences of integers $l_n^-$ and $l_n^+$ such that $l_n^-
\sim l_n^+\sim \sqrt{n}t^{-1}/\sqrt{2}$ and the probability of the
event $\{l_n^-\leq k_n\leq l_n^+\}$ tends to $1$ as $n\rightarrow
\infty$. But on this event, we have by monotonicity
\begin{equation*}
  f\left(\frac{1}{N}\mathcal{X}_{n-l_n^-}'^{[N]}\right)\leq
  f\left(\frac{1}{N}\mathcal{X}_{n-k_n}'^{[N]}\right)\leq 
  f\left(\frac{1}{N}\mathcal{X}_{n-l_n^+}'^{[N]}\right).  
\end{equation*}
The expectations of the outer quantities converge to the right side
of~\eqref{eq:onedimconv-2}. This finishes the proof.
\end{prooof}

\subsection{Convergence of ladder epochs}
\label{rwconvergence}
Aldous and Pitman have shown in~\cite{ALDPIT} that the exponential time
change 
\begin{equation*}
F(t) = \mathfrak{X}(-\ln t),\quad t>0, 
\end{equation*}
with $F(0) = (1,0,\ldots)$ transforms the standard additive coalescent into
a fragmentation process which is self-similar with index $\alpha
=1/2$. In~\cite{BER:2}, Bertoin has given an explicit construction of this
fragmentation process in terms of ladder epochs of Brownian excursion with
drift, and our result on finite-dimensional convergence for the ternary
coalescent will be based on this identity.
 
Let us introduce some notation. We denote by $C[0,1]$ the space of
continuous real-valued paths on $[0,1]$, endowed with the uniform
topology. For an arbitrary path $\omega\in C[0,1]$, its ladder
time set is given by
\begin{equation*}
\mathcal{L}(\omega) = \left\{s\in [0,1]: \omega(s) =
\inf_{[0,s]}\omega\right\}.
\end{equation*}
Since $\mathcal{L}(\omega)$ is a closed set, there exists a unique
decomposition of $[0,1]\backslash\mathcal{L}(\omega)$ into a countable
union of disjoint (open) intervals. We denote by $G(\omega)$ the ranked
sequence of their lengths. By filling up with zeros, we may always
interpret $G(\omega)$ as a mass partition in $\massp$. Note that
$G(\omega)\in\mathbb{S}_1$ if and only if $\mathcal{L}(\omega)$ has
Lebesgue measure zero. 

The construction of the dual fragmentation process
$F$ in~\cite{BER:2} can be summarized as follows. Let $\epsilon =
(\epsilon(s),\, 0\leq s\leq 1)$ be a positive Brownian excursion. For every
$t\geq 0$, consider the excursion dragged down with drift $t$, that is
$\epsilon_t(s)=\epsilon(s)-st$, $0\leq s\leq 1$, and its ladder time set
$\mathcal{L}(\epsilon_t)$, which has almost surely Lebesgue measure
zero. Then, the law of $(G(\epsilon_t),\, t\geq 0)$ and
$(F(t),\,t\geq 0)$ coincide.

In light of our representation of the ternary coalescent in terms of ladder
epochs, it seems natural to establish convergence of these objects. In this
direction, the main step is to prove convergence of the underlying random
paths, with the origin placed at the first instant when their minimum is
attained, towards a Brownian excursion with drift.

To begin with, take a process $(J_n(t),\, t\geq 0)$ distributed as
$(J^{[N]}(t/N^{3/2}),\,t\geq 0)$, and independently of this a Markov chain
$\{X_l\}_{0\leq l\leq n}$ as defined in Section~\ref{randomevolution}. Let
us first fix $t>0$, and write $J_n=J_n(t)$. Remember that given $J_n$, we may
identify $X_{J_n}$ with simple random walk up to time $N$, conditioned to
end at $-(2(n-J_n)+1)$,
\begin{equation*}
S(X_{J_n})_j = 2\left(\sum_{i=0}^{j-1}X_{J_n}(i)\right) - j,\quad
0\leq j\leq 2n+1.
\end{equation*}
By linear interpolation, we define the
correspon\-ding continuous random path $S_{n,t}$ on the unit interval,
\begin{equation*}
S_{n,t}(s) = 2\left(\sum_{i=0}^{\lfloor
  Ns\rfloor-1}X_{J_n}(i)+(Ns-\lfloor Ns\rfloor)X_{J_n}(\lfloor
  Ns\rfloor)\right) - Ns,\quad 0\leq s\leq 1.
\end{equation*}
We shall now prove convergence of the finite-dimensional laws
of the $C[0,1]$-valued process $(N^{-1/2}S_{n,t},\, t>0)$. The
limiting object $(B_{t^{-1}}^{br},\, t>0)$ is distributed as
\begin{equation}
\label{eq:limobj}
 (B_{t^{-1}}^{br},\,t>0) \eql \left((B^{br}(s) -st^{-1},\, 0\leq s\leq 1),\, t> 0\right),
\end{equation}
where $B^{br}$ is a standard Brownian bridge on the unit interval. In
particular, for each fixed $t$, the distribution of $B_{t^{-1}}^{br}$ on
$C[0,1]$ is that of a Brownian bridge from $0$ to $-t^{-1}$.
\begin{lemma}
The $C[0,1]$-valued process $\left(N^{-1/2}S_{n,t},\,t>0\right)$
converges in the sense of finite-dimensional distributions as
$n\rightarrow\infty$ to $\left(B_{t^{-1}}^{br},\, t>0\right)$.
\end{lemma}
\begin{prooof}
  Let us fix $t>0$ as above and first prove one-dimensional convergence. For $0\leq
  s\leq 1$, define
\begin{equation*}
W_n(s) = 2\left(\sum_{i=0}^{\lfloor Ns\rfloor-1}X_{n}(i)+ (Ns-\lfloor
  Ns\rfloor)X_{n}(\lfloor
  Ns\rfloor)\right) - Ns,
\end{equation*}
\begin{equation*}
D_n(s) = 2\left(\sum_{i=0}^{\lfloor
  Ns\rfloor-1}(X_n(i)-X_{J_n}(i)) +(Ns-\lfloor Ns\rfloor)\left(X_{n}(\lfloor Ns\rfloor)-X_{J_n}(\lfloor Ns\rfloor)\right)\right).
\end{equation*}
We may then express $S_{n,t}$ as $S_{n,t} = W_n -D_n$.

The process $W_n(\cdot)$ is linear interpolation of simple random walk up
to time $N$, conditioned to end at $-1$. We deduce from a conditioned
version of Donsker's invariance principle (see Dwass and
Karlin~\cite{DWKA}) that $(N^{-1/2}W_n(s),\, 0\leq s\leq 1)$ converges
weakly in $C[0,1]$ to the standard Brownian bridge $B^{br}$.

Concerning the drift part $D_n$, we let
\begin{eqnarray*}
D^{(1)}_n(s) &=& \sum_{i=0}^{\lfloor
  Ns\rfloor-1}\left(X_n(i)-X_{J_n}(i)\right),\\
 D^{(2)}_n(s) &=&
2(Ns-\lfloor Ns\rfloor)\left(X_{n}(\lfloor Ns\rfloor)-X_{J_n}(\lfloor Ns\rfloor)\right),
\end{eqnarray*}
so that $D_n = 2D_n^{(1)} + D_n^{(2)}$. Now fix $s\in [0,1]$. A moment's
thought reveals that conditioned on $J_n =n-k$ for some $k\in\{0,\ldots, n\}$, the random
variable $D^{(1)}_n(s)$ follows the hypergeometric distribution. More
precisely,
\begin{equation*}
  \pP\left(D^{(1)}_n(s) = j\left|\right.\,J_n=n-k\right) = \frac{{\lfloor Ns\rfloor\choose
      j}{N-\lfloor Ns\rfloor\choose k-j}}{{N\choose k}},
\end{equation*}
where $\max\{0,
  k+\lfloor Ns\rfloor -N\}\leq j\leq\min\{k,\lfloor Ns\rfloor\}$.
As a consequence, 
\begin{equation}
\label{eq:hypergeom-EandV}
\pE\left[D^{(1)}_n(s)\left|\right.\,J_n=n-k\right] = k\frac{\lfloor
  Ns\rfloor}{N},\quad\mbox{Var}\left(D^{(1)}_n(s)\left|\right.\,J_n=n-k\right) \leq k.
\end{equation}
Let $k_n=n-J_n$. Choosing $\varepsilon >0$ arbitrarily small, we have for
large $n$ by the law of total probability
\begin{eqnarray*}
\lefteqn{\pP\left(N^{-1/2}|D_n(s)-2k_ns| > \varepsilon\right)}\\ 
&\leq& \sum_{k=0}^{\lfloor \sqrt{n}t^{-1}\rfloor}\pP\left(N^{-1/2}|D_n^{(1)}(s)-\pE[D^{(1)}_n(s)| > \varepsilon/3
\left|\right.\,k_n=k\right)\pP(k_n=k)\\
&&\ +\ \pP\left(k_n
\geq\sqrt{n}t^{-1}\right)\\
&=&o(1),
\end{eqnarray*}
where the last line follows
from~\eqref{eq:onedimconv-1},~\eqref{eq:hypergeom-EandV} and Chebyshev's
inequality.
Since by~\eqref{eq:asymplemma-1}, $N^{-1/2}2k_n s$ converges in probability
to $t^{-1}s$, so does $N^{-1/2}D_n(s)$. In particular, the
finite-dimensional laws of $(N^{-1/2}D_n(s),\, 0\leq s\leq 1)$
converge to those of $(t^{-1}s,\, 0\leq s\leq 1)$.
Moreover, $D_n(s)$ is increasing in $s$, and a similar computation entails
that for $\lambda$ large enough, as $n\rightarrow \infty$,
\begin{equation*}
\pP\left(N^{-1/2}D_n(1) \geq \lambda\right) = o(1).
\end{equation*}
By Theorem 8.4 of Billingsley~\cite{BILL}, we conclude that the
distributions of $N^{-1/2}D_n(\cdot)$ form a tight sequence. It
follows that $(N^{-1/2}D_n(s),\, 0\leq s\leq 1)$ converges in probability to
$(t^{-1}s,\, 0\leq s\leq 1)$. Applying now Theorem 4.4 from~\cite{BILL}
together with the continuous mapping theorem finishes the proof of the
one-dimensional convergence.

The arguments obviously extend to finite-dimensional distributions. Indeed,
the bridge term $W_n$ is the same for all $t$, and the drift term $D_n$
converges in probability, for each $t$. Therefore, finite-dimensional
convergence follows again from Theorem 4.4 of~\cite{BILL}.
\end{prooof}
Similar as for discrete paths, we introduce for $v\in[0,1]$ the shift
operator $\theta$ on $C[0,1]$,
\begin{equation*}
 (\theta_v\omega)(s) = \left\{\begin{array}{l@{\ ,\ \ }l}
     \omega(s+v)-\omega(v)  & 0\leq s\leq 1-v\\
     \omega(s+v-1) - \omega(v) + \omega(1)-\omega(0)&
     1-v<s\leq 1\end{array}\right..
\end{equation*}
Define $H : C[0,1]\rightarrow[0,1]$ as the first time when the global
minimum is attained,
\begin{equation*}
H(\omega) = \inf\left\{s\in[0,1] : \omega(s) =\inf_{[0,1]}\omega\right\}.
\end{equation*}
Clearly, $H$ is not continuous on the whole space, but it is so restricted
to the subset of paths which uniquely attain their minimum. It is
well-known and also implied by the subsequent Lemma~\ref{localminima-unique}
that the distribution of $B_{t^{-1}}^{br}$ is fully supported on
this subset.
Further, the shift operator is continuous as a map $\theta: C[0,1]\times
[0,1]\rightarrow C[0,1]$, $\theta(\omega,v) = \theta_v\omega$. Setting
$\theta_H\omega = \theta_{H(\omega)}\omega$, it then follows from the
above lemma and the continuous mapping theorem that for
$n\rightarrow\infty$,
\begin{equation*}
\left(N^{-1/2}\theta_HS_{n,t},\,t>0\right) \rightarrow
\left(\theta_HB_{t^{-1}}^{br},\,t>0\right)
\end{equation*} 
in the sense of finite-dimensional distributions.
Recall $(B_{t^{-1}}^{br}(s),\,0\leq s\leq 1)\eql (B^{br}(s) -st^{-1},\, 0\leq
s\leq 1)$, where $B^{br}$ is a Brownian bridge (the same for all
$t$).  Denoting by $\epsilon$ a standard Brownian excursion, it has been
proven by Vervaat in~\cite{VER} that
\begin{equation*}
\theta_HB^{br} \eql \epsilon.
\end{equation*}
Since $\theta_u\circ\theta_v = \theta_w$ for $w=u+v[\textup{mod }1]$, we have
$\theta_H = \theta_H\circ\theta_v$ pathwise for every $0\leq v\leq 1$.
Therefore, if $\mu$ denotes the almost surely unique instant when $B^{br}$
attains its minimum,
\begin{eqnarray}
\label{eq:bridge-excursion-vervaat}
\theta_HB_{t^{-1}}^{br} &\eql&
\theta_H\circ\theta_{\mu}\left(B^{br}-st^{-1},\, 0\leq s\leq 1\right)\nonumber\\
&=& \theta_H\left(\theta_{H}B^{br}-st^{-1},\, 0\leq s\leq 1\right)\nonumber\\
&\eql&\theta_H\epsilon_{t^{-1}}.
\end{eqnarray}
Here, as above, $\epsilon_{t^{-1}}(s) = \epsilon(s)
-st^{-1}$ is the Brownian excursion dragged down with drift $t^{-1}$. Since
$\epsilon_{t^{-1}}$ attains its minimal value almost surely at the endpoint,
we have proven the following
\begin{corollary}
\label{rwpathconv}
In the notation above, $\left(N^{-1/2}\theta_HS_{n,t},\,t>0\right)$
converges in the sense of finite-dimensional distributions as
$n\rightarrow\infty$ to $\left(\epsilon_{t^{-1}},\, t>0\right)$.
\end{corollary}
The convergence of the ternary coalescent is now easy to establish. As
last preparation, let us recall a technical result.
Call a point $x\in [0,1]$ a {\it local
  minimum} of $\omega\in C[0,1]$, if there exists $\delta >0$ such that for
all $y\in [\max\{x-\delta,0\},\min\{x+\delta,1\}]$, $\omega(x)\leq
\omega(y)$. The following statement is true for all real $t$.
\begin{lemma}
\label{localminima-unique}
With probability one, all local minima of $(\varepsilon_t(s),\, 0\leq s\leq
1)$ are distinct.
\end{lemma}
\begin{prooof}
By~\eqref{eq:bridge-excursion-vervaat}, we may show the statement for
$(B_{t}^{br}(s),\,0\leq s\leq 1)$ instead. Since for the time-reversed
process, it holds that
\begin{equation*}
(B^{br}(1-s) -(1-s)t,\, 0\leq
s\leq 1) \eql (B^{br}(s) +st - t,\, 0\leq
s\leq 1),
\end{equation*}
it suffices to show that for some $1/2\leq r<1$, $(B_{t}^{br}(s),\,0\leq s\leq r)$ has almost
surely distinct local minima. However, if $\pF_r$ denotes the filtration
generated by the canonical process $x_{\cdot}$ on $C[0,1]$ up to time
$r<1$, $\pQ$ denotes the 
law of $B_{t}^{br}$ and, for a moment, $\pP$ is Wiener measure and $p$ the
Gaussian transition kernel, it is well-known that $\pQ$ is locally absolute
continuous with respect to $\pP$,
\begin{equation*}
\pQ|_{\pF_r}= \frac{p_{1-r}(x_r,-t)}{p_1(0,-t)}\cdot \pP|_{\pF_r}.
\end{equation*}
Since the local minima of Brownian motion on $[0,1]$ are distinct almost
surely (see for example Theorem 2.11 in the book of M\"orters and
Peres~\cite{MP}), the lemma is proven.
\end{prooof}

\begin{prooof2}{\bf of Theorem~\ref{convtoaddcoal-thm}:}
In view of Bertoin's result in~\cite{BER:2},
the claim follows if we show that the finite-dimensional laws of
\begin{equation*}
t\mapsto \frac{1}{N}\mathcal{X}^{[N]}(t/N^{3/2}),\quad t > 0,
\end{equation*}
converge to those of $(G(\varepsilon_{t^{-1}}),\, t > 0)$. Remember the map
$\vp$ constructed in Section~\ref{mc-and-fcts} sending configurations to
mass partitions. With $J_n(t)=J^{[N]}(t/N^{3/2})$ defined as above, we
have already seen that
\begin{equation*}
\left(\frac{1}{N}\vp(X_{J_n(t)}),\, t\geq 0\right) \eql
  \left(\frac{1}{N}\mathcal{X}^{[N]}(t/N^{3/2}),\, t \geq 0\right). 
\end{equation*}
Let $t>0$, and assume that conditionally
on $J_n$,
\begin{equation*}
\frac{1}{N}\vp(X_{J_n(t)}) =
(s_1,\ldots,s_{2(n-J_n(t))+1}),
\end{equation*}
where $Ns_i\in\{1,3,5,\ldots,N\}$ with $\sum s_i =1$.
Then by construction of both $\vp$, $G$ and linear interpolation,
\begin{equation*}
G(N^{-1/2}\theta_HS_{n,t}) =
(g_1,\ldots,g_{2(n-J_n(t))+1}), 
\end{equation*}
with $g_i=s_i-1/N$ for all $i$.
Thus, the theorem follows if we show finite-dimensional convergence of
$(G(N^{-1/2}\theta_HS_{n,t}),\, t>0)$ to $(G(\varepsilon_{t^{-1}}),\,
t > 0)$. It is easy to check that $G:C[0,1]\rightarrow \massp$
is continuous on the subset of those paths which attain 
their local minima at unique points.
By Lemma~\ref{localminima-unique}, the distribution of
$\varepsilon_{t^{-1}}$ assigns mass one to this subset. Therefore,
Corollary~\ref{rwpathconv} and the continuous mapping theorem yield
convergence of the finite-dimensional distributions on $\massp$, and since
$G(\varepsilon_{t^{-1}})\in\mathbb{S}_1$ with probability one, we
obtain finite-dimensional convergence on $\mathbb{S}_1$.
\end{prooof2}

\subsubsection{Concluding remarks}
(i) For the proof of Theorem~\ref{convtoaddcoal-thm} we used the random walk
representation.  Let us point out another possibility to derive
convergence, using the random binary forest representation.  Following the
construction in Section~\ref{binaryforests}, the state chain of the ternary
coalescent starting from $N$ particles of unit mass can be realized in
reversed time by deleting successively pairs of outgoing edges from a
random tree uniformly distributed over all binary plane trees on $N$
vertices. Such a random tree can be seen as a Galton-Watson tree with
offspring distribution $\mu(k) = \frac{1}{2}(\delta_0(k) +\delta_2(k))$,
conditioned to have total population size $N$. One finds oneself in the
setting of Theorem 23 (in the sublattice case) of Aldous~\cite{ALD3}. In
particular, if $\tau^{[N]}$ denotes the uniform binary plane tree on $N$
vertices, where mass $1/N$ is assigned to each vertex and the edges are
rescaled to have length $1/\sqrt{N}$, then $\tau^{[N]}$ converges weakly as
$N\rightarrow \infty$ to the Brownian continuum random tree (CRT)
introduced in~\cite{ALD1}. By splitting the skeleton of this tree into
subtrees according to a Poisson process of cuts with some intensity $t\geq
0$ per unit length, Aldous and Pitman~\cite{ALDPIT} derived from the CRT an
$\mathbb{S}_1$-valued fragmentation process of ranked masses of tree
components, indexed by the intensity $t$. Further, they showed that the
time change $t\mapsto e^{-t}$ turns this process into the standard additive
coalescent.  Similar as in~\cite{ALDPIT}, it should be possible to
approximate the Poisson process of marks on the CRT by the process of
deleting edges from the binary plane tree. This would lead to another proof
of Theorem~\ref{convtoaddcoal-thm}.\\
(ii) Recall the random walk representation introduced in
Section~\ref{section-coalfragproc}. Fix an integer $k$ of size at least
$3$, and define the configuration space $\cs_n^k$
as the set of all subsets of $\{0,\ldots,(k-1)n\}$ with
cardinality less or equal to $n$. Now identify a configuration $x\in\cs_n^{k}$
with a path of a walk that goes up $k-2$ steps if a site is
occupied and one step down otherwise, i.e.
$S^{(k)}(x)_0 = 0$ and for $1\leq j\leq (k-1)n+1$,
\begin{equation*}
S^{(k)}(x)_j =  
k\left(\sum_{i=0}^{j-1}x(i)\right) - j.
\end{equation*}

By imposing an analogous dynamics, i.e. by occupying successively $n$ sites
chosen uniformly at random from $\{0,\ldots,(k-1)n\}$, the sequence of
ladder epochs of the corresponding new paths is now a realization of the
state chain of the $k$-ary coalescent process with kernel
$\kappa_k(r_1,\ldots,r_k)=r_1+\ldots+r_k +k/(k-2)$, started from $(k-1)n+1$
particles of unit mass. As for the case $k=3$, running this process
backwards in time yields a fragmentation process. Moreover, Kemperman's
formula applies also to first hitting times of such asymmetric random
walks, so that their distributions can easily be computed. With some minor
modifications, and under a different rescaling of time, one again obtains
convergence of the finite-dimensional laws of this $k$-ary coalescent
process towards those of the standard additive coalescent.

Not surprisingly, there is an analogous random $(k-1)$-ary forest
representation of this process. Indeed, when glueing (full) $(k-1)$-ary
trees by picking uniformly at random one leaf and $k-1$ roots from
different components, in a similar way as described in
Section~\ref{binaryforests} for the case $k=3$, the ranked sequence of the
tree sizes is another realization of the state chain of the $k$-ary
coalescent with kernel $\kappa_k$.

This remark shows that our ternary coalescent process is only one particular process out of a family of $k$-ary coalescents
that can be studied by the same means.

\subsection*{Acknowledgments}
I am grateful to Jean Bertoin for introducing me to the topic and helpful
advice. Further I would like to thank two anonymous referees for their
valuable comments.




\begin{thebibliography}{Lam00}
\bibitem{ALD1} Aldous, D. J.: {\it The Continuum Random Tree I.}
  Ann. Probab. {\bf 19}(1) (1991), 1-28. 
\bibitem{ALD3} Aldous, D. J.: {\it The Continuum Random Tree III.}
 Ann. Probab. {\bf 21}(1) (1993), 248-289.
\bibitem{ALDPIT} Aldous, D. J., Pitman, J.: {\it The standard additive
  coalescent.} Ann. Probab. {\bf 26}(4) (1998), 1703-1726. 
\bibitem{BER:2} Bertoin, J.: {\it A fragmentation process connected to
  Brownian motion.} Probab. Theory Related Fields {\bf 117}(2) (2000),
  289-301.
\bibitem{BER:4} Bertoin, J.: {\it Eternal additive coalescents and certain
    bridges with exchangeable increments.} Ann. Probab. {\bf 29}(1) (2001), 344-360. 
\bibitem{BER:3} Bertoin, J.: {\it Different aspects of a random
  fragmentation model.} Stochastic Processes and their Applications {\bf 116}(3) (2006), 345-369.
\bibitem{BER:1} Bertoin, J.: {\it Random Fragmentation and Coagulation
    Processes.} Cambridge Studies in Advanced Mathematics, No. {\bf 102} (2006).
\bibitem{BCP} Bertoin, J., Chaumont, L., Pitman, J.: {\it Path
    transformations of first passage bridges.}
  Elect. Comm. in Probab. {\bf 8} (2003), 155-166.

\bibitem{BILL} Billingsley, P.: {\it Convergence of Probability Measures.}
  Wiley-Interscience (1968).
\bibitem{BRE} Breiman, L.: {\it Probability.} Addison-Wesley, Reading,
  Massachusetts (1968).
\bibitem{DWKA} Dwass, M., Karlin, S.: {\it Conditioned limit theorems.}
  Ann. Math. Stat. {\bf 34}(4) (1963), 1147-1167 .
\bibitem{EVPIT} Evans, S. N., Pitman, J.: {\it Construction of Markovian
  coalescents.} Ann. Inst. H. Poincar\'e, Probab. Statist. {\bf 34}(3) (1998), 339-383.
\bibitem{GK} Gnedenko, B. V., Kolmogorov, A. N.: {\it Limit Distributions for
    Sums of Independent Random Variables.} Addison-Wesley, Reading,
  Massachusetts (1968).
\bibitem{KMP} Kemperman, J. H. B.: {\it The passage problem for a stationary Markov chain.} Statistical Research
Monographs, Vol. {I}, The University of Chicago Press, Chicago (1961).
\bibitem{LG} Le Gall, J.-F.: {\it Random trees and applications.}
  Probability Surveys {\bf 2} (2004), 245-311.
\bibitem{MP} M\"orters, P., Peres, Y.: {\it Brownian motion.} Cambridge
  University Press (2010).
\bibitem{PPY} Perman, M., Pitman, J., Yor, M.: {\it Size biased sampling of
    Poisson point processes and excursions.}
   Probab. Theory Related Fields {\bf 92}(1) (1992), 21-39.
\bibitem{PIT:99} Pitman, J.: {\it Coalescent Random Forests.}
  J. Combin. Theory Ser. A {\bf 85}(2) (1999), 165-193.
\bibitem{PIT:StF} Pitman, J.: {\it Combinatorial Stochastic
    Processes. \'Ecole d'\'et\'e de Probabilit\'es de St. Flour.} Lecture
  Notes in Mathematics {\bf 1875}, Springer (2006).
\bibitem{VER} Vervaat, W.: {\it A relation between Brownian bridge and
  Brownian excursion.} Ann. Probab. {\bf 7}(1) (1979), 141-149.
\end{thebibliography}
\end{document}